\documentclass[twocolumn]{autart}    
\usepackage{graphicx}
\usepackage{amsmath,amssymb,mathtools,graphicx}
\usepackage{arydshln}
\usepackage{listings,color,caption}
\usepackage{amsmath}
\usepackage{bm}
\usepackage{cite}
\usepackage{verbatim}

\usepackage{etoolbox}
\newbool{arxivversion}
\setbool{arxivversion}{true}

\usepackage{algpseudocode}
\usepackage{algorithm}

\def\JV{\textcolor{black}}

\begin{document}

\begin{frontmatter}
%\runtitle{Insert a suggested running title}  % Running title for regular 
                                              % papers but only if the title  
                                              % is over 5 words. Running title 
                                              % is not shown in output.

\title{Robust targeted exploration for systems with non-stochastic disturbances\thanksref{footnoteinfo}} % Title, preferably not more 
                                                % than 10 words.

\thanks[footnoteinfo]{Frank Allg\"ower is thankful that his work was funded by Deutsche Forschungsgemeinschaft (DFG, German Research Foundation) under Germany’s Excellence Strategy - EXC 2075 - 390740016 and under grant 468094890. Frank Allg\"ower acknowledges the support by the Stuttgart Center for Simulation Science (SimTech). Janani Venkatasubramanian thanks the International Max Planck Research School for Intelligent Systems (IMPRS-IS) for supporting her. Johannes K\"ohler was supported by the Swiss National Science Foundation under NCCR Automation (grant agreement 51NF40 180545).}

\author[1]{Janani Venkatasubramanian},
%\ead{janani.venkatasubramanian@ist.uni-stuttgart.de},
\author[2]{Johannes K\"ohler},
% \ead{jkoehle@ethz.ch},
\author[3]{Mark Cannon},
% \ead{mark.cannon@eng.ox.ac.uk},
\author[1]{Frank Allg\"ower}
% \ead{frank.allgower@ist.uni-stuttgart.de}

\address[1]{Institute for Systems Theory and Automatic Control, University of Stuttgart, 70550 Stuttgart, Germany; $\quad\quad\quad\quad\quad$ (e-mail: janani.venkatasubramanian, frank.allgower@ist.uni-stuttgart.de)}  % Please supply                             
\address[2]{Institute for Dynamic Systems and Control, ETH Z\"urich, Z\"urich CH-80092, Switzerland; (e-mail: jkoehle@ethz.ch)}             % full addresses
\address[3]{Control Group, University of Oxford, Parks
Road, Oxford OX1 3PJ, United Kingdom; (e-mail: mark.cannon@eng.ox.ac.uk)}        % here.

                             % keyword list or with the 
                                          % help of the Automatica 
                                          % keyword wizard

\begin{abstract}
We propose a novel targeted exploration strategy designed specifically for uncertain linear time-invariant systems with energy-bounded disturbances, i.e., without any assumptions on the distribution of the disturbances. We use classical results characterizing the set of non-falsified parameters consistent with energy-bounded disturbances. We derive a semidefinite program which computes an exploration strategy that guarantees a desired accuracy of the parameter estimate. This design is based on sufficient conditions on the spectral content of the exploration data that robustly account for initial parametric uncertainty. Finally, we highlight the applicability of the exploration strategy through a numerical example involving a nonlinear system.
% The proposed exploration strategy guarantees a desired error-bound on the estimated parameters after finite-time exploration.
% This work utilizes classical data-dependent uncertainty bounds on the least-squares parameter estimates in the presence of energy-bounded noise.
% Specifically, we derive sufficient conditions on the spectral content of the exploration data that ensures a desired error bound on the estimated parameters. 

% We utilize these sufficient conditions to derive robust LMIs, while accounting for initial parametric uncertainty. 

% This leads to a semidefinite program-based design which computes the exploration inputs that guarantee a desired error-bound on the parameters. 
\end{abstract}

\end{frontmatter}

\section{Introduction}\label{sec:intro}
Designing reliable controllers for unknown dynamical systems requires accurate knowledge of the model parameters, which can be obtained from data \cite{gevers2005identification}. The accuracy of the parameters significantly depends on the quality of the data used for system identification. Informative data can be strategically obtained from an experiment through the process of \textit{targeted exploration} or \textit{optimal experiment design}~\cite{pronzato2008optimal, gevers1986optimal}. Specifically, targeted exploration inputs are tailored to reduce model uncertainty, thereby ensuring the attainment of a desired accuracy in the identified model \cite{jansson2005input, bombois2021robust}, or the feasibility of robust control design ~\cite{barenthin2008identification, umenberger2019robust, ferizbegovic2019learning, venkatasubramanian2023sequential}. In this paper, we develop a targeted exploration strategy for uncertain linear systems subject to energy-bounded disturbances.

% Classical system identification offers asymptotic methods that produce accurate models, but they require an infinite amount of data to achieve this level of precision~\cite{ljung1999system}. However, in practice, we only have finite data from experiments, and the reliability of the model-based controller depends significantly on the quality of the data used for system identification. In order to identify accurate models, informative data can be strategically obtained from a finite experiment through the process of \textit{targeted exploration} or \textit{optimal experiment design}~\cite{pronzato2008optimal, gevers1986optimal}. Specifically, targeted exploration inputs are tailored to reduce model uncertainty, thereby ensuring the attainment of a desired accuracy in the identified model \cite{jansson2005input, bombois2021robust}, or the feasibility of robust control design ~\cite{barenthin2008identification, umenberger2019robust, ferizbegovic2019learning, venkatasubramanian2020robust, venkatasubramanian2023sequential}. In this paper, we develop a finite targeted exploration strategy for uncertain linear systems with energy-bounded disturbances.

Targeted exploration is typically studied in a stochastic setup, for example, with independent and identically distributed (i.i.d.) disturbances with zero mean. In this case, one can construct a data-based confidence ellipsoid for the parameters \cite{ljung1999system} that can be approximately predicted and optimized. This classical asymptotic result has been utilized to design targeted exploration strategies for dual control methods \cite{barenthin2008identification, umenberger2019robust,iannelli2020structured,venkatasubramanian2023sequential}. In~\cite{dean2017sample}, tight confidence regions for the parameters are constructed, assuming the data is independent, and hence it is not applicable to correlated time-series data from a single trajectory. 
Nevertheless, this bound has been utilized in the design of targeted exploration for dual control in \cite{ferizbegovic2019learning}.
% While these approaches offer robust performance guarantees, they rely on asymptotic identification results, i.e., models identified from infinitely long experiments. In contrast to these asymptotic results, there has been a recent resurgence in the non-asymptotic (finite-sample) analysis of system identification \cite{abbasi2011regret, dean2017sample}, supported by advances in statistical learning theory \cite{tsiamis2023statistical}. Analogous to \cite{umenberger2019robust}, finite-sample data-based high probability credibility regions for parameter estimates are introduced in \cite{dean2017sample}. In particular, the uncertainty bounds on the model parameters in \cite{dean2017sample} require data to be independent, and hence not applicable to correlated time-series data. Nevertheless, these bounds have been utilized in the design of targeted exploration for dual control in \cite{ferizbegovic2019learning}.
% In contrast to these asymptotic results, there has been a resurgence in the study of finite-sample system identification, specifically to improve robust and adaptive linear quadratic regulator (LQR) design since the recent initiation by \cite{abbasi2011regret, dean2017sample}. 
%This result stems from the `coarse-ID control' procedure in which a model and its associated uncertainty is estimated from a few experimental trials, and is then utilized for controller design.5
A common feature among all of the discussed targeted exploration approaches is that they consider linear systems subject to i.i.d. zero mean stochastic disturbances. However, real-world systems often exhibit nonlinear behaviour or unmodeled dynamics which introduce additional deterministic model mismatch and cannot be explained by independent stochastic noise~\cite{sarker2023accurate}. Instead, such effects can be modeled as bounded disturbances, assuming that the disturbances belong to a known bounded set, e.g., point-wise bounds or energy bounds. To quantify models under such disturbances, set-membership estimation methods have been developed, and various identification results in this direction can be found in \cite{fogel1979system, milanese1989estimation, sarker2023accurate, bisoffi2021trade, shakouri2024identification}. Although these methods have recently gained popularity for data-driven robust control design \cite{vanwaarde2022noisy, berberich2022combining}, a principled method for optimal experiment design in the presence of bounded disturbances is still lacking.

% Frequently, the statistics of the disturbance are either unknown or inaccurately assumed, leading to potentially invalid probabilistic credibility regions.
% \textbf{Related work:} The majority of relevant literature focuses on linear time-invariant (LTI) systems and deal with linear system identification and learning for the linear quadratic regulator (LQR).  

In this paper, we design a targeted exploration strategy that ensures a desired error bound on the estimated parameters by utilizing a data-dependent uncertainty bound based on energy-bounded disturbances \cite{fogel1979system}. Unlike existing experiment design methods, which often rely on stochastic disturbance assumptions \cite{umenberger2019robust, ferizbegovic2019learning, venkatasubramanian2023sequential}, our approach provides robustness against non-stochastic, adversarial disturbances and is applicable to systems with bounded nonlinearities. We consider multi-sine exploration inputs of specific frequencies and optimized amplitudes to explicitly shape and reduce uncertainty in a targeted manner. As a main contribution, we derive sufficient conditions on the spectral content of the exploration data that asymptotically guarantees a desired error-bound on the parameters estimated through exploration. We utilize these sufficient conditions to derive LMIs for exploration which ensure the desired error bound on the parameters. This approach gives rise to a targeted exploration design with minimal input energy based on a semidefinite program (SDP).
\JV{Unlike existing methods for stochastic disturbances \cite{venkatasubramanian2023sequential,umenberger2019robust,ferizbegovic2019learning}, the non-stochastic identification results considered in this work yield structurally different requirements on excitation, which require new worst-case robust guarantees rather than high-probability guarantees.} To this end, we formulate the problem in Section \ref{sec:probstatement}, and review data-dependent uncertainty bounds in Section \ref{sec:prelim}. In Section \ref{sec:exploration}, we derive the exploration strategy by leveraging a sufficient conditions on time-series exploration data, as proposed in our previous work \cite{venkatasubramanian2023towards}, while additionally robustly accounting for parametric uncertainty and error due to disturbances with suitable bounds. In Section \ref{sec:numerical}, we provide a numerical example to highlight the applicability of the proposed exploration strategy to nonlinear systems. \JV{We remark that the proofs of Lemma \ref{lemma:ub_Zlb_tight}, Propositions \ref{prop2} and \ref{prop4}, and Theorem \ref{thm:finalthm} are not included in this brief paper, but can be found in the technical report \cite{venkatasubramanian2024robust}.}

\section{Problem statement}\label{sec:probstatement}
\textbf{Notation:} For a matrix $A \in \mathbb{R}^{n \times m}$, $A^\top$ denotes the transpose. For a matrix $A \in \mathbb{C}^{n \times m}$, $A^\mathsf{H}$ denotes the conjugate transpose. The positive (semi-)definiteness of a matrix $A \in \mathbb{C}^{n \times n}$ is denoted by ($A = A^\mathsf{H} \succeq 0$) $A = A^\mathsf{H} \succ 0$. The operator $\mathrm{vec}(A)$ stacks the columns of $A$ to form a vector. The operator $\textnormal{diag}(A_1,\dots,A_n)$ creates a block diagonal matrix by aligning the matrices $A_1,\dots,A_n$ along the diagonal with $A_1$ in the upper left corner. The Kronecker product is denoted by $\otimes$. The Euclidean norm and the weighted Euclidean norm for a vector $x\in\mathbb{R}^n$ and a matrix $P \succ 0$ are denoted by $\|x\|=\sqrt{x^\top x}$ and $\|x\|_P=\sqrt{x^\top P x}$, respectively. The largest singular value of a matrix $A\in\mathbb{C}^{m\times n}$ is $\|A\|$. Furthermore, given a matrix $M \succeq 0$, $\|A\|_M = \| M^{1/2} A \|$ where $M^{1/2}$ is the symmetric square root matrix of $M$. The identity matrix of size $n$ is denoted by $I_n$. A vector of ones of size $n$ is denoted by $\mathbf{1}_n \in \mathbb{R}^{n \times 1}$.

\textit{Amplitude of a spectral line \cite[Def. 2]{venkatasubramanian2023sequential}:} The amplitude of the spectral line of a sequence $\{\phi_k\}_{k=0}^{T-1}$ at a frequency $\omega \in \Omega_T:=\{0,1/T,\dots,(T-1)/T\}$ is given by
\begin{equation}\label{eq:spectral_amplitude}
\bar{\phi}(\omega):=\frac{1}{T}\sum_{k=0}^{T-1}\phi_k e^{-j2\pi \omega k}.
\end{equation}
% Furthermore, for a matrix $A\in\mathbb{C}^{m\times n}$, given a matrix $M \in\mathbb{C}^{m\times m}$ and $ M \succeq 0$, $\|A\|_M = \| M^{1/2} A \|$ where $M^{1/2}$ is the unique Hermitian square-root of $M$.
% Given a sequence $\{ x_k\}_{k=0}^{T-1}$, the discrete Fourier transform (DFT) of the sequence is denoted by $\mathbf{x}(e^{j\omega})=\sum_{k=0}^{T-1}x_k e^{-j2\pi k \omega}$ where $\omega \in \Omega_T:=\{0,1/T,\dots,(T-1)/T\}$. The inverse discrete Fourier transform is denoted by $x_k=\frac{1}{T} \sum_{i=1}^{T} \mathbf{x}(e^{j\omega_i})e^{j2\pi \omega_i k}$ where $\omega_i \in \Omega_T$.

% \textit{Amplitude of a spectral line \cite[Def. 2]{sarker2022parameter}:} Given a sequence $\{\phi_k\}_{k=0}^{T-1}$, the amplitude of the spectral line of the sequence $\bar{\phi}(\omega_0)$ at a frequency $\omega_0 \in \Omega_T:=\{0,1/T,\dots,(T-1)/T\}$ is given by
% $\bar{\phi}(\omega_0):=\frac{1}{T}\sum_{k=0}^{T-1}\phi_k e^{-j2\pi \omega_0 k}$.
% A random variable $X \in \mathbb{R}$ is said to be sub-Gaussian \cite{rigollet2015high}, with variance proxy $\sigma^2 \in \mathbb{R}$, i.e., $X \sim \mathsf{subG}(\sigma^2)$, if $\mathbb{E}[X]=0$ and its moment generating function satisfies
% \begin{equation}
% \mathbb{E}\left[e^{sX}\right] \leq e^{\left( \frac{\sigma^2 s^2}{2} \right)},\quad \forall s \in \mathbb{R}.
% \end{equation}
\subsection{Setting}
Consider a discrete-time, linear time-invariant system of the form
\begin{align}\label{eq:sys}
x_{k+1}=A_\mathrm{tr} x_k + B_\mathrm{tr} u_k + w_k,
\end{align}
 where $x_k \in \mathbb{R}^{n_\mathrm{x}}$ is the state, $u_k \in \mathbb{R}^{n_\mathrm{u}}$ is the control input, and $w_k \in \mathbb{R}^{n_\mathrm{x}}$ is the disturbance. 
 In our setting, the true system parameters $A_\mathrm{tr}$, $B_\mathrm{tr}$ are uncertain. Hence, it is necessary to collect informative data from an optimal experiment for a fixed $T\in \mathbb{N}$ time steps to improve parameter accuracy. \JV{To simplify the exposition, we assume that the initial state is at the origin, i.e., $x_0=0$.} Furthermore, we assume that the state can be measured and the disturbances are energy-bounded.
 \begin{assum}\label{a0} The disturbances $w$ are energy-bounded; there exists a known constant $\gamma_\mathrm{w}> 0$ such that 
 \begin{align}\label{eq:noisebound}
     \sum_{k=0}^{T-1}\|w_k\|^2 \leq \gamma_\mathrm{w}.
\end{align}\vspace{-2em}
\end{assum}
This assumption allows the system to exhibit nonlinear behaviour, as long as the nonlinearities can be bounded in energy (as discussed in Section \ref{sec:numerical}). Importantly, this assumption distinguishes our approach from existing experiment design methods which are typically restricted to i.i.d. Gaussian disturbances, e.g., \cite{umenberger2019robust, ferizbegovic2019learning, venkatasubramanian2023sequential}.

\textbf{Exploration goal:} Since the true system parameters $\theta_\mathrm{tr} =\text{vec}([A_\mathrm{tr}$, $B_\mathrm{tr}]) \in \mathbb{R}^{n_\theta}$, with $n_\theta = n_\mathrm{x}(n_\mathrm{x}+n_\mathrm{u})$, are not precisely known, exploratory inputs should be designed to excite the system to gather informative data. Specifically, our objective is to design inputs that excite the system in a manner as to obtain an estimate $\hat{\theta}_T=\text{vec}([\hat{A}_T,\hat{B}_T])$ that satisfies 
\begin{align}\label{eq:exp_goal2}
    (\theta_\mathrm{tr}-\hat{\theta}_T)^\top (D_\mathrm{des}\otimes I_{n_\mathrm{x}}) (\theta_\mathrm{tr}-\hat{\theta}_T)\leq 1,
\end{align}
where $D_\mathrm{des}\succ 0$ is a user-defined matrix characterizing closeness of $\hat{\theta}_T$ to $\theta_\mathrm{tr}$. We assume that we have some prior knowledge about the system dynamics.
 \begin{assum} \label{a1} The unknown parameters $\theta_\mathrm{tr}=\mathrm{vec}([A_\mathrm{tr},B_\mathrm{tr}])$ lie in a known set $\mathbf{\Theta}_0$, i.e., $\theta_\mathrm{tr} \in \mathbf{\Theta}_0$, where
\begin{align}\label{eq:Theta0}
\mathbf{\Theta}_0:=\left\{\theta:(\hat{\theta}_0-\theta)^\top (D_0 \otimes I_{n_\mathrm{x}})  (\hat{\theta}_0-\theta) \leq 1 \right\},
\end{align}
with an estimate $\hat{\theta}_0=\mathrm{vec}([\hat{A}_0,\hat{B}_0])$ for some $D_0\succ 0$.
\end{assum}
\begin{rem}
An initial estimate $\hat{\theta}_0$ and set $\mathbf{\Theta}_0$ as in Assumption \ref{a1} result naturally from a finite-horizon experiment under Assumption \ref{a0} (cf. Lemma \ref{lem:thetat}).
\end{rem}
% In what follows, we first introduce the form of the exploration inputs.
% Given $\hat{\theta}_T$ and $D_\mathrm{des} \succ 0$, we can define an ellipsoid $\bm{\hat{\Theta}}_T$ centered around the estimate $\hat{\theta}_T$ such that $\theta_\mathrm{tr} \in \bm{\hat{\Theta}}_T$.

\subsection{Exploration strategy}
The exploration input sequence takes the form
\begin{align}\label{eq:exploration_controller}
u_k=\sum_{i=1}^{L} \bar{u}(\omega_i) \cos(2 \pi \omega_i k), \quad\, k=0,\dots,T-1
\end{align}
where $T$ is the exploration time and $\bar{u}(\omega_i) \in \mathbb{R}^{n_\mathrm{u}}$ are the amplitudes of the sinusoidal inputs at $L \in \mathbb{N}$ distinct selected frequencies $\omega_i \in \Omega_T$ with $n_\mathrm{x}+n_\mathrm{u}\leq L \leq T$. In practice, frequencies may be selected based on prior information \cite{mehra1974optimal}. Since the input signal is deterministic and sinusoidal, the amplitude of the spectral line of the sequence $\{u_k\}_{k=0}^{T-1}$ at frequency $\omega_i$ is $\bar{u}(\omega_i)$.
%We select $n_\mathrm{\phi}$ frequencies from $\Omega_T$ and denote $\Omega_\phi$ as the set of all $n_\mathrm{\phi}$ selected frequencies.
Denote $U_\mathrm{e}=\mathrm{diag}(\bar{u}(\omega_1),\dots,\bar{u}(\omega_L)) \in \mathbb{R}^{Ln_\mathrm{u} \times L}$. The exploration input is computed such that it excites the system sufficiently with minimal input energy, based on the initial parameter estimates. Bounding the input energy by a constant $\gamma^2_\mathrm{e}$ is equivalent to $\sum_{i=1}^L \|\bar{u}(\omega_i)\|^2= \mathbf{1}_{L}^\top U_\mathrm{e}^\top U_\mathrm{e} \mathbf{1}_{L} \preceq \gamma_\mathrm{e}^2$ where $\mathbf{1}_{L} \in \mathbb{R}^{L\times 1}$ is a vector of ones, and the bound $\gamma_\mathrm{e} \geq 0$ is desired to be small. Using the Schur complement, this criterion is equivalent to 
\begin{equation}\label{eq:min_energy_cost}
S_{\textnormal{energy-bound}}(\gamma_\mathrm{e},U_\mathrm{e})\coloneqq \begin{bmatrix} 
\gamma_\mathrm{e} & \mathbf{1}_{L}^\top U_\mathrm{e}^\top \\ U_\mathrm{e} \mathbf{1}_{L} & \gamma_\mathrm{e} I
\end{bmatrix} \succeq 0.
\end{equation}
In order to design the exploration inputs, we make the following assumption regarding the system dynamics. 
 \begin{assum}\label{a2}The system matrix $A_\mathrm{tr}$ is Schur stable.
\end{assum}
\begin{rem} \label{rem:2} $A_\mathrm{tr}$ is required to be Schur stable since we consider only open-loop inputs in our exploration strategy \eqref{eq:exploration_controller}. Assumption \ref{a2} could be relaxed by augmenting the inputs in \eqref{eq:exploration_controller} with an additional linear feedback, i.e., $v_k=u_k+Kx_k$, ensuring robust stability for all $\theta \in \mathbf{\Theta}_0$ \eqref{eq:Theta0}.
\end{rem}
To achieve the exploration goal, the amplitudes of the sinusoidal exploration inputs need to be optimized such that by applying the exploration inputs, the obtained estimate satisfies the desired uncertainty bound~\eqref{eq:exp_goal2}. 

% To this end, we first provide preliminaries regarding uncertainty bounds for parameter estimation based on time-series data, and spectral information in Section \ref{sec:prelim}. As one of the main results, we derive sufficient conditions in terms of spectral information of the exploration data, and the targeted exploration strategy in Section \ref{sec:exploration}. Consequently, upon solving the exploration problem and applying the exploration inputs, we identify parameter estimates up to the desired accuracy relative to the true parameters \eqref{eq:exp_goal2}.

% \subsubsection*{Exploration goal} The goal of the exploration strategy is to design exploration inputs that excite the system in a manner as to obtain an estimate $\hat{\theta}_T=\mathrm{vec}([\hat{A}_T,\hat{B}_T])$ that satisfies
% \begin{align}\label{eq:exp_goal2}
%     (\theta_\mathrm{tr}-\hat{\theta}_T)^\top D_\mathrm{des} (\theta_\mathrm{tr}-\hat{\theta}_T)\leq 1,
% \end{align}
% where $D_\mathrm{des} \succ 0$ is a user-defined matrix characterizing the closeness of $\hat{\theta}_T$ to $\theta_\mathrm{tr}$. Given $\hat{\theta}_T$ and $D_T \succ 0$, we can define a credibility region $\bm{\hat{\Theta}}_T$ centered around the estimate $\hat{\theta}_T$ such that $\theta_\mathrm{tr} \in \bm{\hat{\Theta}}_T$ where
% \begin{align}
%     \bm{\hat{\Theta}}_T:=\left\{\theta:(\theta-\hat{\theta}_T)^\top D_\mathrm{des} (\theta-\hat{\theta}_T)\leq 1\right\}.
% \end{align}

\section{Preliminaries on data-driven uncertainty quantification}\label{sec:prelim}
% % Furthermore, we discuss the derivation of frequency domain information from time-series data that can be used to compute the uncertainty bound on the unknown parameters.

% \subsection{Uncertainty bound}\label{sec:uncertaintybound_prelim}
In this section, we discuss a data-dependent uncertainty bound on the parameter estimates in the presence of energy-bounded noise \cite{fogel1979system}. Given observed data $\mathcal{D}_{T+1}=\{x_k,u_k\}_{k=0}^{T}$ of length $T+1$, the objective is to quantify the uncertainty associated with the unknown parameters $\theta_\mathrm{tr}$. Henceforth, we denote $\phi_k=[x_k^\top\, u_k^\top]^\top \in \mathbb{R}^{n_\phi}$ where $n_\phi=n_\mathrm{x}+n_\mathrm{u}$. The system \eqref{eq:sys} can be re-written in terms of parameter $\theta_\mathrm{tr}=\textnormal{vec}([A_\mathrm{tr},B_\mathrm{tr}])$ as
 \begin{align*}
 x_{k+1}=(\phi_k^\top \otimes I_{n_\mathrm{x}})\theta_\mathrm{tr} + w_k.
 \end{align*}
In order to simplify the exposition, we denote
\begin{align*}
    \Phi&=[\phi_0,\dots,\phi_{T-1}] \in \mathbb{R}^{n_\phi\times T},\\
    X^\top& =[x_1^\top,\dots, x_T^\top] \in \mathbb{R}^{1 \times Tn_\mathrm{x}}.
\end{align*}
We obtain the following expressions for the mean $\hat{\theta}_T~=~\textnormal{vec}([\hat{A}_T,\hat{B}_T])$ and covariance $P$ of the parameters from the standard least squares formulation \cite[Section 1.3]{ljung1999system}:
\begin{align}\label{eq:mean_est}
    \hat{\theta}_T&=P\sum_{k=0}^{T-1} (\phi_k^\top \otimes I_{n_\mathrm{x}})^\top x_{k+1}=P (\Phi\otimes I_{n_\mathrm{x}}) X,\\
    \label{eq:covar_est}
    P&=\left( \sum_{k=0}^{T-1} \phi_k \phi_k^\top \right)^{-1} \otimes I_{n_\mathrm{x}}=(\Phi \Phi^\top)^{-1} \otimes I_{n_\mathrm{x}}.
\end{align}

The non-falsified region for the uncertain parameters $\theta$ is provided in the following lemma.

\begin{lem}\label{lem:thetat} \cite{fogel1979system} Let Assumption~\ref{a0} hold. Given data set $\mathcal{D}_{T+1}$, the set of non-falsified parameters $\theta$ is given by
\begin{align}\label{eq:Theta_T}
    \mathbf{\Theta}_T:=\left\{ \theta: (\theta-\hat{\theta}_T )^\top{P^{-1}} (\theta-\hat{\theta}_T) \leq G \right\}
\end{align}
where
% \begin{align}\label{eq:G}
%     G =\gamma_\mathrm{w} + \|\hat{\theta}_T\|_{P^{-1}}^2-\sum_{k=0}^{T-1}\|x_{k+1}\|^2.
% \end{align}
\begin{align}\label{eq:G}
    G =\gamma_\mathrm{w} + \|\hat{\theta}_T\|_{P^{-1}}^2-X^\top X.
\end{align}
\end{lem}
\ifbool{arxivversion}{\begin{pf}
The energy constraint on the process noise in \eqref{eq:sys} yields the following non-falsified set:
\begin{align}
\mathbf{\Theta}_T&=\left\{ \theta: \sum_{k=0}^{T-1} \|x_{k+1}-(\phi_k^\top \otimes I_{n_\mathrm{x}})\theta \|^2 \leq \gamma_\mathrm{w} \right\}
% &=\left\{ \theta: (\theta-\hat{\theta}_T )^\top{P_{\theta}^{-1}} (\theta-\hat{\theta}_T) \leq G(T)\right\}.
\end{align}
which can be equivalently written as
\begin{align}\label{eq:Gtbound_pre}
\nonumber
    &\theta^\top \left((\Phi \Phi^\top) \otimes I_{n_\mathrm{x}}  \right)\theta 
    -2\left(X^\top (\Phi^\top \otimes I_{n_\mathrm{x}})  \right)\theta\\ \leq & \;\gamma_\mathrm{w}- X^\top X.
\end{align}   
By adding $\|\hat{\theta}_T\|_{P^{-1}}^2$ to both sides of \eqref{eq:Gtbound_pre}, and using \eqref{eq:mean_est} and \eqref{eq:covar_est} to complete the squares, we get
\begin{align}\label{eq:Gtbound}
    \|\theta-\hat{\theta}_T\|_{P^{-1}}^2 \leq \gamma_\mathrm{w} + \|\hat{\theta}_T\|_{P^{-1}}^2-X^\top X =:G,
\end{align}
which is equivalent to \eqref{eq:Theta_T}. $\hfill\square$
\end{pf}}{}
Under Assumption \ref{a0}, the non-falsified set $\bm{\Theta}_T$ provides an exact characterization of the set of parameters explaining the data. The ellipsoid \eqref{eq:Theta_T}, derived from energy-bounded constraints, is characterized by $\hat{\theta}_T$ and $P$, which coincide with the mean and covariance of the least squares estimator under Gaussian disturbances \cite[Prop. 2.1]{umenberger2019robust}. However, unlike the Gaussian case, the scaling $G$ \eqref{eq:G} of the bounding ellipsoid is also data-dependent in the case of energy-bounded disturbances, which significantly affects optimal experiment design. Furthermore, for $\gamma_\mathrm{w} \propto T$, $G$ scales at most linearly with $T$ as $T \to \infty$ since $G \leq \gamma_\mathrm{w}$ \cite[Lemma 4]{fogel1979system}. In contrast, the scaling of the confidence ellipsoid in the Gaussian case does not depend on $T$ \cite{ljung1999system, umenberger2019robust}. Since $P^{-1}$ increases linearly with $T$, the confidence ellipsoid in the Gaussian case shrinks with $T$~\cite{ljung1999system,umenberger2019robust}. However, in the considered case of energy-bounded disturbances, the size of the non-falsified set $\mathbf{\Theta}_T$, in general, does not decrease as $T \to \infty$.

\section{Targeted Exploration}\label{sec:exploration}
In this section, we propose a targeted exploration strategy based on the data-dependent uncertainty bound in Lemma \ref{lem:thetat}. The exploration strategy builds upon sufficient conditions on the exploration data outlined in Section \ref{subsec:suffcon}. In particular, we derive bounds on the exploration data using the spectral information of the exploration inputs in Section \ref{subsec:lowerboundfe}. Since these bounds depend on the uncertain model parameters, bounds on the effect of model uncertainty are derived in Section \ref{subsec:hinf_bound}. We utilize these bounds to derive sufficient conditions on the spectral information of the exploration inputs in Section \ref{sec:suff_spec}. However, since the derived sufficient conditions are non-convex in the decision variables, a convex relaxation procedure is carried out in Section \ref{subsec:convexrel}. Finally, in Section \ref{subsec:finalsuffcon}, the exploration problem is reduced to a set of LMIs that provide us exploration inputs that ensure the exploration goal.

\subsection{Sufficient conditions for exploration}\label{subsec:suffcon}

Given the form of the exploration inputs in \eqref{eq:exploration_controller}, the exploration goal \eqref{eq:exp_goal2}, and the data-dependent uncertainty bound in Lemma~\ref{lem:thetat}, in what follows, we provide conditions that the exploration data have to satisfy to achieve the exploration goal.  Denote the Cholesky decomposition of $D_\mathrm{des}$ as $D_\mathrm{des}=D_\mathrm{des}^{\frac{1}{2}\top}D_\mathrm{des}^{\frac{1}{2}}$. The following proposition presents a sufficient condition to ensure that the exploration goal is achieved.

\begin{thm} \label{thm:sufficientcond1} \JV{\cite[Theorem 4]{venkatasubramanian2023towards}}
Suppose $\Phi$ and $X$ satisfy
\begin{align}\label{eq:pe5}
\nonumber
    \begin{bmatrix}
        \Phi \Phi^\top - \gamma_\mathrm{w} D_\mathrm{des}& 0\\
        0 & 0
    \end{bmatrix}&\\ 
    +\underbrace{\begin{bmatrix}
       D_\mathrm{des}^{\frac{1}{2}\top}(X^\top \otimes I_{n_\phi}) \\ (\Phi \otimes I_{n_\mathrm{x}}) \otimes I_{n_\phi}
    \end{bmatrix}}_{=:Z} \begin{bmatrix}
       D_\mathrm{des}^{\frac{1}{2}\top}(X^\top \otimes I_{n_\phi}) \\ (\Phi \otimes I_{n_\mathrm{x}}) \otimes I_{n_\phi}
    \end{bmatrix}^\top &\succeq 0.
\end{align}
Then, the estimate $\hat{\theta}_T$ computed as in \eqref{eq:mean_est} satisfies the exploration goal \eqref{eq:exp_goal2}.
\end{thm}
\ifbool{arxivversion}{\begin{pf}
The bound in \eqref{eq:Gtbound} can be re-written as
\begin{align}\label{eq:Gtbound_re}
\nonumber
    &(\theta-\hat{\theta}_T)^\top ((\Phi \Phi^\top) \otimes I_{n_\mathrm{x}}) (\theta-\hat{\theta}_T)\\
\nonumber
    \overset{\eqref{eq:mean_est}}{\leq} & \gamma_\mathrm{w}-X^\top( I - (\Phi^\top\otimes I_{n_\mathrm{x}}) P (\Phi\otimes I_{n_\mathrm{x}}))X\\
    \overset{\eqref{eq:covar_est}}{=} & \gamma_\mathrm{w}-X^\top X +X^\top ((\Phi^\top (\Phi \Phi^\top)^{-1}\Phi) \otimes I_{n_\mathrm{x}} )X.
\end{align}
By applying the Schur complement twice to \eqref{eq:Gtbound_re}, we get
\begin{align}\label{eq:bound_schur}
\nonumber
    &(\theta-\hat{\theta}_T) (\theta-\hat{\theta}_T)^\top\\ 
    \nonumber \preceq & (\gamma_\mathrm{w}-X^\top X +X^\top ((\Phi^\top (\Phi \Phi^\top)^{-1}\Phi) \otimes I_{n_\mathrm{x}} )X) \\& \cdot (\Phi \Phi^\top)^{-1}\otimes I_{n_\mathrm{x}}.
\end{align}
Inequality \eqref{eq:pe5} can be written as
\begin{small}
\begin{equation}\label{eq:schur_pe}
    \left[\begin{array}{c;{2pt/2pt}c}
       \begin{matrix}\Phi \Phi^\top-\gamma_\mathrm{w}D_\mathrm{des}\\
        +D_\mathrm{des}^{\frac{1}{2}\top}\left((X^\top X) \otimes I_{n_\phi}\right)D_\mathrm{des}^{\frac{1}{2}}
       \end{matrix}
        & \star^\top \\ \hdashline[2pt/2pt]
        \underbrace{(((\Phi \otimes I_{n_\mathrm{x}}) X)\otimes I_{n_\phi})D_\mathrm{des}^{\frac{1}{2}}}_{\star} & ((\Phi \Phi^\top)\otimes I_{n_\mathrm{x}})\otimes I_{n_\phi}
    \end{array}\right] \succeq 0.
\end{equation}
\end{small}
By applying the Schur complement to \eqref{eq:schur_pe}, we get
\begin{small}
\begin{align}\label{eq:pe24}
 \Phi \Phi^\top -D_\mathrm{des}^{\frac{1}{2}\top}\left((\gamma_\mathrm{w} - X^\top X)\otimes I_{n_\phi}\right)D_\mathrm{des}^{\frac{1}{2}}& \\\nonumber
    -D_\mathrm{des}^{\frac{1}{2}\top}\left((X^\top ((\Phi^\top (\Phi \Phi^\top)^{-1}\Phi) \otimes I_{n_\mathrm{x}}  ) X ) \otimes I_{n_\phi} \right)D_\mathrm{des}^{\frac{1}{2}} &\succeq 0.
\end{align}
\end{small}
Since $\Phi \Phi^\top \succeq 0 \iff (\Phi \Phi^\top)\otimes I_{n_\mathrm{x}} \succeq 0$, \eqref{eq:pe24} can be written as \begin{align}\label{eq:pe_ineq}
\nonumber
    (\Phi \Phi^\top)\otimes I_{n_\mathrm{x}} \succeq & \big(X^\top ((\Phi^\top (\Phi \Phi^\top)^{-1}\Phi) \otimes I_{n_\mathrm{x}} )X \\
    &+\gamma_\mathrm{w}-X^\top X\big) (D_\mathrm{des} \otimes I_{n_\mathrm{x}}).
\end{align}
Furthermore, by inserting \eqref{eq:pe_ineq} in \eqref{eq:bound_schur}, we get 
\begin{align}\label{eq:expgoal_schur}
    (\theta-\hat{\theta}_T) (\theta-\hat{\theta}_T)^\top \preceq (D_\mathrm{des} \otimes I_{n_\mathrm{x}})^{-1}.
\end{align}
Finally, applying the Schur complement twice to \eqref{eq:expgoal_schur} yields the exploration goal \eqref{eq:exp_goal2}. $\hfill\square$
\end{pf}}{}
To compute the amplitudes of the exploration inputs $U_\mathrm{e}$, a few challenges need to be addressed. Note that Inequality \eqref{eq:pe5} depends on $X$ and $\Phi$ quadratically, which further depend on the amplitudes of the exploration inputs $U_\mathrm{e}$ \eqref{eq:exploration_controller} and the disturbance $w$.  Furthermore, due to uncertain true dynamics $A_\mathrm{tr},\,B_\mathrm{tr}$, the linear mapping from the input sequence to the state sequence is not known. 

\subsection{Bounds based on the theory of spectral lines}\label{subsec:lowerboundfe} 
In what follows, we address the aforementioned issues by determining sufficient conditions for targeted exploration in terms of the spectral content of $X$ and $\Phi$ based on the theory of spectral lines. Given $u_k$ as in \eqref{eq:exploration_controller}, $x_k$ has $L$ spectral lines from $0$ to $T-1$ at distinct frequencies $\omega_i \in \Omega_T$, $i=1,\dots,L$ with amplitudes \cite[Lemma 1]{sarker2020parameter}:
\begin{align}\label{eq:xbar}
\nonumber
\bar{x}(\omega_i)=&\underbrace{(e^{j 2 \pi \omega_i}I - A_\mathrm{tr})^{-1}B_\mathrm{tr}}_{=:V_{\mathrm{x},i}} \bar{u}(\omega_i)\\
& +\underbrace{(e^{j 2 \pi \omega_i}I - A_\mathrm{tr})^{-1}}_{=:Y_{\mathrm{x},i}} \bar{w}(\omega_i) + \bar{x}_{\mathrm{err}}(\omega_i).
\end{align}
The transient error in the amplitude of a spectral line $\bar{x}_{\mathrm{err}}(\omega_i)$ decays uniformly (with rate $\frac{1}{\sqrt{T}}$) to $0$ as $T \to \infty$ (cf. Assumption \ref{a2}, \cite[Theorem 2.1]{ljung1999system}). To simplify the exposition, we will assume that it can be neglected \cite{sarker2020parameter}.
\begin{assum}\label{a4} The transient error satisfies $\bar{x}_{err}(\omega_i)=0$ for all $\omega_i \in \Omega_T$
\end{assum}
Note that this assumption holds naturally if we let $T~\to~\infty$. We refer the reader to the work in \cite{venkatasubramanian2025beyondcdc} for a rigorous treatment of the transient error term. 
More compactly, let us define
%Although the transient error term can be bounded using the results in \cite{venkatasubramanian2025beyondcdc}, we focus on the asymptotic setting for a streamlined presentation.
\begin{align}\label{eq:Xbar}
\bar{X}=[\bar{x}(\omega_1)^\top,\dots, \bar{x}(\omega_L)^\top]^\top \in \mathbb{C}^{n_\mathrm{x}L}
\end{align}
which satisfies
\begin{align}\label{eq:xbarcompact}
\bar{X}=  \underbrace{V_{\mathrm{x},\mathrm{tr}}U_\mathrm{e}\mathbf{1}_L}_{=:\bar{X}_{\mathrm{u}}}+\underbrace{Y_{\mathrm{x},\mathrm{tr}}W\mathbf{1}_L}_{=:\bar{X}_{\mathrm{w}}},  
\end{align}
with 
\begin{align}\label{eq:VxYxW}
\nonumber
V_{\mathrm{x},\mathrm{tr}}&:=\textnormal{diag}(V_{\mathrm{x},1},\cdots,V_{\mathrm{x},L})\in \mathbb{C}^{n_\mathrm{x} L\times n_\mathrm{u} L},\\
\nonumber
Y_{\mathrm{x},\mathrm{tr}}&:=\textnormal{diag}(Y_{\mathrm{x},1},\cdots,Y_{\mathrm{x},L}) \in \mathbb{C}^{n_\mathrm{x}  L\times n_\mathrm{x} L},\\
W&:=\mathrm{diag}(\bar{w}(\omega_1),\dots,\bar{w}(\omega_L))\in \mathbb{C}^{n_\mathrm{x}L \times L}.
\end{align}
Furthermore, $\phi_k$ has $L$ spectral lines from $0$ to $T-1$ at distinct frequencies $\omega_i \in \Omega_T$, $i=1,\dots,L$ with amplitudes
\begin{align}\label{eq:phibar}
\bar{\phi}(\omega_i):=\underbrace{\begin{bmatrix}
V_{\mathrm{x},i} \\ I_{n_\mathrm{u}}
\end{bmatrix}}_{=:V_{\phi,i}} \bar{u}(\omega_i)+\underbrace{\begin{bmatrix}
Y_{\mathrm{x},i} \\ 0
\end{bmatrix}}_{=:Y_{\phi,i}} \bar{w}(\omega_i).
\end{align}
We compactly define
\begin{align}
    \bar{\Phi}=[\bar{\phi}(\omega_1), \dots ,\bar{\phi}(\omega_L)] \in \mathbb{C}^{n_\phi \times L},
\end{align}
which satisfies
\begin{align}\label{eq:phibarcompact}
    \bar{\Phi} = \underbrace{V_\mathrm{\phi,tr}U_\mathrm{e}}_{=:\bar{\Phi}_\mathrm{u}} + \underbrace{Y_\mathrm{\phi,tr} W}_{=:\bar{\Phi}_\mathrm{w}},
\end{align}
with 
\begin{align}\label{eq:VphiYphi}
\nonumber
V_{\phi,\mathrm{tr}}&:=[V_{\phi,1},\cdots,V_{\phi,L}] \in \mathbb{C}^{n_\mathrm{\phi}\times n_\mathrm{u}L},\\
 Y_\mathrm{\phi,tr}&:=[Y_{\phi,i},\cdots,Y_{\phi,L}] \in \mathbb{C}^{n_\mathrm{\phi}\times n_\mathrm{x}L}.
\end{align}
The spectral content of $Z$ \eqref{eq:pe5} is denoted by $\bar{Z} \in \mathbb{C}^{(n_\phi + n_\mathrm{x} n_\phi^2) \times (n_\mathrm{x} n_\phi L)}$:
\begin{align}\label{eq:Zspectral}
\nonumber
\bar{Z}=&\begin{bmatrix}
        D_\mathrm{des}^{\frac{1}{2}\top} ((\bar{X}_\mathrm{u}+\bar{X}_\mathrm{w})^\mathsf{H} \otimes I_{n_\phi})\\ ((\bar{\Phi}_\mathrm{u}+\bar{\Phi}_\mathrm{w}) \otimes I_{n_\mathrm{x}}) \otimes I_{n_\phi}
    \end{bmatrix}\\
    =&\underbrace{\begin{bmatrix}
        D_\mathrm{des}^{\frac{1}{2}\top} ( (\mathbf{1}_L^\top U_\mathrm{e}^\top V_{\mathrm{x,tr}}^\mathsf{H}) \otimes I_{n_\phi})\\  
        ( V_\mathrm{\phi,tr}U_\mathrm{e} ) \otimes I_{n_\mathrm{x}n_\phi}
\end{bmatrix}}_{=:\bar{Z}_\mathrm{u}}+
\underbrace{\begin{bmatrix}
        D_\mathrm{des}^{\frac{1}{2}\top} (\bar{X}_\mathrm{w}^\mathsf{H}  \otimes I_{n_\phi})\\
        \bar{\Phi}_\mathrm{w} \otimes I_{n_\mathrm{x}n_\phi}
    \end{bmatrix}}_{=:\bar{Z}_\mathrm{w}}.
\end{align}
The following lemma provides lower bounds on $\Phi \Phi^\top$ and $ZZ^\top$ \eqref{eq:pe5} using the spectral content of the signals $x_k$ and $\phi_k$.

\begin{lem} \label{lemma:ub_Zlb_tight}
    Let Assumptions \ref{a2} and \ref{a4} hold. For any $\epsilon \in (0,1)$, $\phi_k$ and $Z$ satisfy
\begin{align}\label{eq:phiphitight}
   \Phi \Phi^\top \succeq T \left(  (1-\epsilon) \bar{\Phi}_\mathrm{u} \bar{\Phi}_\mathrm{u}^\mathsf{H} - \left(\frac{1-\epsilon}{\epsilon}\right) \bar{\Phi}_\mathrm{w} \bar{\Phi}_\mathrm{w}^\mathsf{H}  \right),
\end{align}
and
    \begin{align}\label{eq:ZZtight}
    Z Z^\top \succeq T \left((1-\epsilon)\bar{Z}_\mathrm{u} \bar{Z}_\mathrm{u}^\mathsf{H}-\left(\frac{1-\epsilon}{\epsilon} \right)\bar{Z}_\mathrm{w} \bar{Z}_\mathrm{w}^\mathsf{H}  \right),
\end{align}
respectively.
\end{lem}
\ifbool{arxivversion}{The proof of Lemma \ref{lemma:ub_Zlb_tight} is provided in Appendix \ref{app:spectralbounds}.}{\JV{The proof of Lemma \ref{lemma:ub_Zlb_tight} may be derived using the Parseval-Plancheral identity; cf. \cite[Appendix A]{venkatasubramanian2024robust}.}} The matrices $\bar{\Phi}_\mathrm{u}$, $\bar{\Phi}_\mathrm{w}$, $\bar{Z}_\mathrm{u}$ and $\bar{Z}_\mathrm{u}$ in Lemma \ref{lemma:ub_Zlb_tight} depend on the transfer matrices $V_{\mathrm{x,tr}}$, $V_{\mathrm{\phi,tr}}$, $Y_{\mathrm{x,tr}}$, and $Y_{\mathrm{\phi,tr}}$. These transfer matrices are dependent on the true dynamics $A_\mathrm{tr}$, $B_\mathrm{tr}$, and hence, uncertain. Therefore, in what follows, suitable bounds are derived.

\subsection{Bounds on transfer matrices}\label{subsec:hinf_bound}
Denote 
\begin{align}\label{eq:vtilde}
    \tilde{V}_\phi=V_\mathrm{\phi,tr}-\hat{V}_\phi,\; \tilde{V}_\mathrm{x}=V_\mathrm{\mathrm{x},tr}-\hat{V}_\mathrm{x},
\end{align}
where the estimates
\begin{align}\label{eq:VhatVxhat}
\nonumber
\hat{V}_\phi&=[\hat{V}_{\phi,1},\cdots,\hat{V}_{\phi,L}] \in \mathbb{C}^{n_\mathrm{\phi}\times Ln_\mathrm{u}},\\
\hat{V}_\mathrm{x}&=\text{diag}\left(\hat{V}_{\mathrm{x},1},\cdots,\hat{V}_{\mathrm{x},L}\right) \in \mathbb{C}^{n_\mathrm{\mathrm{x}}L\times n_\mathrm{u}L}
\end{align}
are computed using the initial estimates $\hat{\theta}_0=\mathrm{vec}([\hat{A}_0,\hat{B}_0])$ (cf. Assumption \ref{a1}). 
We can compute matrices $\tilde{\Gamma}_\phi$, $\tilde{\Gamma}_\mathrm{x}$, $\Gamma_\mathrm{\phi}$ and $\Gamma_\mathrm{\mathrm{x}}$ such that
\begin{align}\label{eq:tf_prop}
\tilde{V}_\phi\tilde{V}_\phi^\mathsf{H} \preceq \tilde{\Gamma}_{
\phi},&\; \tilde{V}_\mathrm{x} \tilde{V}_\mathrm{x}^\mathsf{H} \preceq \tilde{\Gamma}_{\mathrm{x}},\\
\label{eq:tfy_prop}
Y_\mathrm{\phi,tr}Y_\mathrm{\phi,tr}^\mathsf{H}  \preceq \Gamma_\mathrm{\phi},&\; Y_\mathrm{\mathrm{x},tr}Y_\mathrm{\mathrm{x},tr}^\mathsf{H} \preceq \Gamma_\mathrm{x}
\end{align}
using $\theta_\mathrm{tr} \in \mathbf{\Theta}_0$ (cf. Assumption \ref{a1}). Conditions \eqref{eq:tf_prop} and \eqref{eq:tfy_prop} are LMIs, and hence matrices $\tilde{\Gamma}_\mathrm{x}, \tilde{\Gamma}_\phi \succ 0$ and $\Gamma_\phi,\Gamma_\mathrm{\mathrm{x}} \succ 0$ may be computed using robust control methods as shown in \cite[Appendices A-B]{venkatasubramanian2023sequential}, or scenario methods as shown in  \cite[Appendix C]{venkatasubramanian2023sequential}. Utilizing \eqref{eq:tfy_prop}, we derive bounds on $\bar{\Phi}_\mathrm{w}$, $\bar{X}_\mathrm{w}$ and $\bar{Z}_\mathrm{w}$ in the following lemma.
\begin{lem}\label{lem:noisebounds}
    Let Assumptions \ref{a0}, \ref{a1}, \ref{a2} and \ref{a4} hold. Given the bounds on $Y_{\phi,\mathrm{tr}}$ and $Y_\mathrm{x,tr}$ \eqref{eq:tfy_prop}, we have
    \begin{align}\label{eq:PhiwPhiwH}
        \bar{\Phi}_\mathrm{w} \bar{\Phi}_\mathrm{w}^\mathsf{H} & \preceq \bar{W}_\phi:= \frac{\gamma_\mathrm{w}}{T} \Gamma_\phi
    \end{align}
and
\begin{align}\label{eq:barWZ}
    \bar{Z}_\mathrm{w} \bar{Z}_\mathrm{w}^\mathsf{H} & \preceq \underbrace{\left(\frac{\gamma_\mathrm{w}}{T}\|\Gamma_\mathrm{x}\|\|D_\mathrm{des}\| +     \frac{\gamma_\mathrm{w}}{T}\|\Gamma_\phi\| \right)I_{(n_\phi + n_\mathrm{x} n_\phi^2)}}_{\bar{W}_Z:=}.
\end{align} 
\end{lem}
\ifbool{arxivversion}{The proof of Lemma \ref{lem:noisebounds} is provided in Appendix \ref{app:noisebounds}.}{\JV{Lemma \ref{lem:noisebounds} uses properties of the induced matrix norm; cf. \cite[Appendix B]{venkatasubramanian2024robust}.}} Given the bounds in Lemma \ref{lem:noisebounds}, we next derive spectral lines-based sufficient conditions for targeted exploration.

\subsection{Sufficient conditions for targeted exploration based on the theory of spectral lines}\label{sec:suff_spec}
The following proposition provides a condition in terms of the spectral content of $\phi$, which, if satisfied, ensures that the exploration goal \eqref{eq:exp_goal2} is achieved.
\begin{prop}\label{prop2} Let Assumptions \ref{a0} and \ref{a1} hold. Suppose the matrices $\bar{\Phi}_\mathrm{u}$ and $\bar{Z}_\mathrm{u}$ satisfy
\begin{align}\label{eq:prop2}
\nonumber
    \begin{bmatrix}
        T\left( (1-\epsilon) \bar{\Phi}_\mathrm{u} \bar{\Phi}_\mathrm{u}^\mathsf{H} - \left(\frac{1-\epsilon}{\epsilon}\right) \bar{W}_\phi \right) - \gamma_\mathrm{w} D_\mathrm{des}& 0\\0 & 0
    \end{bmatrix} & \\
    +T\left((1-\epsilon) \bar{Z}_\mathrm{u}\bar{Z}_\mathrm{u}^\mathsf{H} - \left( \frac{1-\epsilon}{\epsilon}\right) \bar{W}_Z \right) & \succeq 0.
\end{align}
Then, the estimate $\hat{\theta}_T$ computed as in \eqref{eq:mean_est} satisfies the exploration goal \eqref{eq:exp_goal2}.
\end{prop}
\ifbool{arxivversion}{The proof of Proposition \ref{prop2} is provided in Appendix \ref{app:prop2}.}{\JV{The proof of Proposition \ref{prop2} can be derived by substituting the terms in Theorem \ref{thm:sufficientcond1} with bounds derived in Lemmas \ref{lemma:ub_Zlb_tight} and \ref{lem:noisebounds}; cf. \cite[Appendix C]{venkatasubramanian2024robust}.}}
Note that $\bar{\Phi}_\mathrm{u}$ and $\bar{Z}_\mathrm{u}$ depend linearly on the decision variable $U_\mathrm{e}$. Determining a lower bound based on the Inequality \eqref{eq:prop2} results in non-convex constraints in $U_\mathrm{e}$. To overcome this problem, we utilize a convex relaxation procedure.

\subsection{Convex relaxation}\label{subsec:convexrel}
%The following lemma is utilized to make Inequality \eqref{eq:prop2} linear in the decision variable $U_\mathrm{e}$.
%\begin{lem}\label{lem:convexrel1} For any matrices $M \in \mathbb{C}^{n \times m}$ and $N \in \mathbb{C}^{n \times m}$, we have
%\begin{align}\label{eq:convex_rel}
%    M M^\mathsf{H} \succeq M N^\mathsf{H} + N M^\mathsf{H} - N N^\mathsf{H}.
%\end{align}
%\end{lem}
%\begin{pf} We have
%$M M^\mathsf{H} - M N^\mathsf{H} - N M^\mathsf{H} + N N^\mathsf{H}=(M-N) (M-N)^\mathsf{H} \succeq 0$
%and hence, \eqref{eq:convex_rel} holds. $\hfill\square$
%\end{pf}

The following proposition provides a sufficient condition linear in $\bar{Z}_\mathrm{u}$ which, ensures the exploration goal \eqref{eq:exp_goal2}.
\begin{prop}\label{prop4} Let Assumptions \ref{a0} and \ref{a1} hold. Suppose the matrices $\bar{\Phi}_\mathrm{u}$, $\bar{Z}_\mathrm{u}$ and $\hat{Z} \in\mathbb{C}^{(n_\phi+n_\mathrm{x} n_\phi^2) \times L n_\mathrm{x} n_\phi}$ satisfy
% $\bar{D} \succeq 0$,
\begin{align}\label{eq:prop3}
\nonumber
    \begin{bmatrix}
        \left( (1-\epsilon) \bar{\Phi}_\mathrm{u} \bar{\Phi}_\mathrm{u}^\mathsf{H} - \left(\frac{1-\epsilon}{\epsilon}\right) \bar{W}_\phi\right)  - \frac{\gamma_\mathrm{w}}{T} D_\mathrm{des}& 0\\0 & 0
    \end{bmatrix} & \\ 
    % -\bar{D} &
    +\left((1-\epsilon) \left( \bar{Z}_\mathrm{u}\hat{Z}^\mathsf{H}+ \hat{Z} \bar{Z}_\mathrm{u}^\mathsf{H} -\hat{Z} \hat{Z}^\mathsf{H}\right)-\left( \tfrac{1-\epsilon}{\epsilon}\right) \bar{W}_Z \right)  & \succeq 0,
\end{align}
then an estimate $\hat{\theta}_T$ computed as in \eqref{eq:mean_est} satisfies the exploration goal \eqref{eq:exp_goal2}.
\end{prop}
\ifbool{arxivversion}{The proof of the proposition is provided in Appendix \ref{app:convexrel}.}{\JV{The proof of the proposition can be derived by substituting the quadratic term $\bar{Z}_\mathrm{u}\bar{Z}_\mathrm{u}^\mathsf{H}$ in Proposition \ref{prop2} with a lower bound linear in $\bar{Z}_\mathrm{u}$ obtained by completing the square; cf. \cite[Appendix D]{venkatasubramanian2024robust}.}} The bound derived in \eqref{eq:prop3} is tight if $\hat{Z}=\bar{Z}_\mathrm{u}$. Since $\bar{Z}_\mathrm{u}$ comprises uncertain elements $V_\mathrm{x,tr}$, $V_\mathrm{\phi,tr}$ and the unknown decision variable $U_\mathrm{e}$, we consider a candidate $\hat{Z}$. Later, this relaxation is embedded an iterative process to reduce conservatism. In what follows, we utilize Inequality \eqref{eq:prop3} to derive a condition linear in the decision variable $U_\mathrm{e}$ that ensures the exploration goal \eqref{eq:exp_goal2}.
 
\subsection{Exploration SDP}\label{subsec:finalsuffcon}
In this section, we provide a sufficient condition that ensure the exploration goal \eqref{eq:exp_goal2} using Proposition \ref{prop4}. In \eqref{eq:Zspectral}, $\bar{Z}_\mathrm{u}$ can be be written as $\bar{Z}_\mathrm{u}=\bar{Z}_\mathrm{u,1}+\bar{Z}_\mathrm{u,2}$, where
\begin{align}
\nonumber
    \bar{Z}_\mathrm{u,1}&=\begin{bmatrix}
        D_\mathrm{des}^{\frac{1}{2}\top}\left( \mathbf{1}_L^\top U_\mathrm{e}^\top \otimes I_{n_\phi} \right)\\
        0
    \end{bmatrix} \left( V_{\mathrm{x,tr}}^\mathsf{H} \otimes I_{n_\phi}\right),\\
    \bar{Z}_\mathrm{u,2}&=\begin{bmatrix}
        0\\ \left(V_\mathrm{\phi,tr}\right)\otimes I_{n_\mathrm{x}n_\phi}
    \end{bmatrix} \left(U_\mathrm{e}\otimes I_{n_\mathrm{x}n_\phi}\right).    
\end{align}
In order to robustly account for uncertainties in Inequality \eqref{eq:prop3}, we split Inequality \eqref{eq:prop3} into three inequalities that are handled separately. The following inequalities are equivalent to Inequality \eqref{eq:prop3} if $\bar{D}_1+\bar{D}_2+\bar{D}_3 \succeq 0$:
 \begin{small}
\begin{subequations}\label{eq:Ds}
\begin{align}
\label{eq:D1}
     (1-\epsilon)\left( \bar{Z}_\mathrm{u,1} \hat{Z}^\mathsf{H} + \hat{Z}\bar{Z}_\mathrm{u,1}^\mathsf{H}\right) -\bar{D}_1&\succeq 0,\\
\label{eq:D2}
     (1-\epsilon)\left( \bar{Z}_\mathrm{u,2} \hat{Z}^\mathsf{H} + \hat{Z}\bar{Z}_\mathrm{u,2}^\mathsf{H} - \hat{Z} \hat{Z}^\mathsf{H}\right)-\left( \tfrac{1-\epsilon}{\epsilon}\right) \bar{W}_Z -\bar{D}_2&\succeq 0,\\
\label{eq:D3}
     \begin{bmatrix}
        \left( (1-\epsilon) \bar{\Phi}_\mathrm{u} \bar{\Phi}_\mathrm{u}^\mathsf{H} - \left(\frac{1-\epsilon}{\epsilon}\right) \bar{W}_\phi \right) - \frac{\gamma_\mathrm{w}}{T} D_\mathrm{des}& 0\\0 & 0
    \end{bmatrix}- \bar{D}_3 & \succeq 0.
 \end{align}   
 \end{subequations} 
 \end{small}
\begin{figure*}
\begin{subequations}
\begin{align}\label{eq:sprocD1}
&S_\textnormal{exp-1}(\epsilon, \tau_1, U_\mathrm{e},\hat{Z}, \hat{V}_\mathrm{x}, \tilde{\Gamma}_\mathrm{x},D_\mathrm{des},\bar{D}_1):=&\\
\nonumber    &\begin{bmatrix}
        0 & \begin{bmatrix}
        (1-\epsilon)D_\mathrm{des}^{\frac{1}{2}\top}\left( \mathbf{1}_L^\top U_\mathrm{e}^\top \otimes I_{n_\phi} \right)\\
        0
    \end{bmatrix}^\mathsf{H} \\
        \begin{bmatrix}
        (1-\epsilon)D_\mathrm{des}^{\frac{1}{2}\top}\left( \mathbf{1}_L^\top U_\mathrm{e}^\top \otimes I_{n_\phi} \right)\\
        0
    \end{bmatrix}& -\bar{D}_1  \end{bmatrix} - \tau_1 \begin{bmatrix}
        -I & (\hat{V}_\mathrm{x}^\mathsf{H}  \otimes I_{n_\phi})\hat{Z}^\mathsf{H}\\
     \hat{Z}(\hat{V}_\mathrm{x} \otimes I_{n_\phi}) & \hat{Z}(( \tilde{\Gamma}_\mathrm{x}-\hat{V}_\mathrm{x} \hat{V}_\mathrm{x}^\mathsf{H})\otimes I_{n_\phi})\hat{Z}^\mathsf{H}
    \end{bmatrix} \succeq 0
\end{align}
\begin{align}\label{eq:sprocD2}
&S_\textnormal{exp-2}(\epsilon, \tau_2, U_\mathrm{e},\hat{Z},\bar{W}_Z,\hat{V}_\phi, \tilde{\Gamma}_\phi,D_\mathrm{des},\bar{D}_2):=&\\
\nonumber
    &\begin{bmatrix}
        0 & (1-\epsilon)(U_\mathrm{e}\otimes I_{n_\mathrm{x} n_\phi}) \hat{Z}^\mathsf{H}\\
        (1-\epsilon)\hat{Z} (U_\mathrm{e}^\top\otimes I_{n_\mathrm{x} n_\phi}) & -(1-\epsilon)\hat{Z} \hat{Z}^\mathsf{H}-\left(\tfrac{1-\epsilon}{\epsilon}\right)\bar{W}_Z - \bar{D}_2
    \end{bmatrix}-\tau_2 \begin{bmatrix}
        -I & \begin{bmatrix}
        0\\ \hat{V}_\mathrm{\phi}\otimes I_{n_\mathrm{x} n_\phi}
    \end{bmatrix}^\mathsf{H}\\
    \begin{bmatrix}
        0\\ \hat{V}_\mathrm{\phi}\otimes I_{n_\mathrm{x} n_\phi}
    \end{bmatrix} & \begin{bmatrix}
        0 & 0\\ 0 & (\tilde{\Gamma}_\phi - \hat{V}_\mathrm{\phi}\hat{V}_\mathrm{\phi}^\mathsf{H})\otimes I_{n_\mathrm{x} n_\phi}
    \end{bmatrix}
    \end{bmatrix}\succeq 0
\end{align}
\begin{align}\label{eq:sprocD3}
&S_\textnormal{exp-3}(\epsilon, \tau_3, U_\mathrm{e},\hat{U},\bar{W}_\phi,\hat{V}_\phi, \tilde{\Gamma}_\phi,\gamma_\mathrm{w}, D_\mathrm{des},\bar{D}_3):=&\\
    \nonumber
& \begin{bmatrix} (1-\epsilon)(U_\mathrm{e} \hat{U}^\top + \hat{U} U_\mathrm{e}^\top - \hat{U} \hat{U}^\top) & 0\\
        0 & \begin{bmatrix}
            -\left(\tfrac{1-\epsilon}{\epsilon}\right)\bar{W}_\phi-\frac{\gamma_\mathrm{w}}{T} D_\mathrm{des}  & 0\\
            0 & 0
        \end{bmatrix} -\bar{D}_{3}
    \end{bmatrix}-\tau_3 \begin{bmatrix}
         -I & \begin{bmatrix}
            \hat{V}_\mathrm{\phi}\\0
        \end{bmatrix}^\mathsf{H}  \\
        \begin{bmatrix}
            \hat{V}_\mathrm{\phi}\\0
        \end{bmatrix}
     &  \begin{bmatrix} \tilde{\Gamma}_\phi - \hat{V}_\mathrm{\phi}\hat{V}_\mathrm{\phi}^\mathsf{H} & 0 \\ 0 & 0 \end{bmatrix}
    \end{bmatrix}\succeq 0
\end{align}
\end{subequations}    
\end{figure*}
The following theorem provides a sufficient condition linear in $U_\mathrm{e}$ which ensures the exploration goal \eqref{eq:exp_goal2}.
\begin{thm}\label{thm:finalthm}
    Let Assumptions \ref{a0}, \ref{a1}, \ref{a2} and \ref{a4} hold. Suppose there exist matrices $U_\mathrm{e}$, $\bar{D}_1$, $\bar{D}_2$ and $\bar{D}_3$, and scalars $\tau_1 \geq 0$, $\tau_2 \geq 0$, and $\tau_3 \geq 0$ such that 
    \begin{align}
    \nonumber
        S_\textnormal{exp-1}(\epsilon, \tau_1, U_\mathrm{e},\hat{Z}, \hat{V}_\mathrm{x}, \tilde{\Gamma}_\mathrm{x},D_\mathrm{des},\bar{D}_1) \succeq 0,\\
    \nonumber    
        S_\textnormal{exp-2}(\epsilon, \tau_2, U_\mathrm{e},\hat{Z},\bar{W}_Z,\hat{V}_\phi, \tilde{\Gamma}_\phi,D_\mathrm{des},\bar{D}_2) \succeq 0,\\
    \nonumber
        S_\textnormal{exp-3}(\epsilon, \tau_3, U_\mathrm{e},\hat{U},\bar{W}_\phi,\hat{V}_\phi, \tilde{\Gamma}_\phi,\gamma_\mathrm{w}, D_\mathrm{des},\bar{D}_3) \succeq 0,\\
        \label{eq:thm12}
    \bar{D}_1+\bar{D}_2+\bar{D}_3 \succeq 0,
    \end{align}
    where $S_\textnormal{exp-1}$, $S_\textnormal{exp-2}$, and $S_\textnormal{exp-3}$ are defined in \eqref{eq:sprocD1}, \eqref{eq:sprocD2} and \eqref{eq:sprocD3}, respectively. Then, an estimate $\hat{\theta}_T$ computed as in \eqref{eq:mean_est} upon the application of the input \eqref{eq:exploration_controller} satisfies the exploration goal \eqref{eq:exp_goal2}.
\end{thm}
\ifbool{arxivversion}{The proof of Theorem \ref{thm:finalthm} is provided in Appendix \ref{app:finalthm}.}{\JV{The proof of Theorem \ref{thm:finalthm} is provided in \cite[Appendix E]{venkatasubramanian2024robust}.}} The key idea of the proof is the application of robust control tools, in particular, the matrix S-lemma~\cite{boyd2004convex,vanwaarde2022noisy}, to account for parametric uncertainty. In particular, feasibility of LMIs (49b) and (49c) requires $(\tilde{\Gamma}_\phi - \hat{V}_\mathrm{\phi}\hat{V}_\mathrm{\phi}^\mathsf{H}) \prec 0$, which holds if the initial uncertainty is sufficiently small. Consequently, we can pose the exploration problem of designing exploration inputs that excite the system with minimal energy to obtain an estimate $\hat{\theta}_T$ that satisfies \eqref{eq:exp_goal2} using the following SDP:
\begin{small}
\begin{align}\label{eq:exp_problem}
\nonumber
    \underset{U_\mathrm{e},\gamma_\mathrm{e}, \atop {\tau_1 \geq 0  , \tau_2 \geq 0, \atop \tau_3 \geq 0}}{\inf}  & \; \gamma_\mathrm{e}\\
\nonumber
    \text{s.t. }& \; S_{\textnormal{energy-bound}}(\gamma_\mathrm{e},U_\mathrm{e})\succeq 0\\
\nonumber
    & \; S_\textnormal{exp-1}(\epsilon, \tau_1, U_\mathrm{e},\hat{Z}, \hat{V}_\mathrm{x}, \tilde{\Gamma}_\mathrm{x},D_\mathrm{des},\bar{D}_1) \succeq 0\\
\nonumber
    & \; S_\textnormal{exp-2}(\epsilon, \tau_2, U_\mathrm{e},\hat{Z},\bar{W}_Z,\hat{V}_\phi, \tilde{\Gamma}_\phi,D_\mathrm{des},\bar{D}_2) \succeq 0\\
\nonumber
    & \; S_\textnormal{exp-3}(\epsilon, \tau_3, U_\mathrm{e},\hat{U},\bar{W}_\phi,\hat{V}_\phi, \tilde{\Gamma}_\phi,\gamma_\mathrm{w}, D_\mathrm{des},\bar{D}_3) \succeq 0\\
    & \; \bar{D}_1+\bar{D}_2+\bar{D}_3 \succeq 0.
\end{align}    
\end{small}
A solution of \eqref{eq:exp_problem} gives us the parameters required for the implementation of the exploration input, i.e., $U_\mathrm{e}=\mathrm{diag}(\bar{u}(\omega_1),\dots,\bar{u}(\omega_L))$, which guarantees the desired uncertainty bound $D_\mathrm{des}$ \eqref{eq:exp_goal2}. The suboptimality introduced by the convex relaxation procedure can be reduced by iterating Problem \eqref{eq:exp_problem} multiple times until $\gamma_\mathrm{e}$ does not change by re-computing $\hat{U}$ and $\hat{Z}$ for the next iteration as
\begin{align}\label{eq:initcand}
    \hat{U}&=U_\mathrm{e}^*,\; \hat{Z}&=\begin{bmatrix}
        D_\mathrm{des}^{\frac{1}{2}\top} \left(\left(\mathbf{1}_L^\top U_\mathrm{e}^{*\top} \hat{V}_{\mathrm{x}}^\mathsf{H} \right)\otimes I_{n_\phi}\right)\\     \left(\hat{V}_\mathrm{\phi}U_\mathrm{e}^*\right)\otimes I_{n_\mathrm{x} n_\phi}
    \end{bmatrix}
\end{align}
wherein $U_\mathrm{e}^*$ is the solution from the previous iteration. 

\begin{rem} (Efficient implementation and conservatism)
	If $D_\mathrm{des}$ is a scaled identity matrix, i.e., $D_\mathrm{des}=cI_{n_\phi}$ with $c>0$, the Kronecker product with $I_{n_\phi}$ in LMIs \eqref{eq:sprocD1} and \eqref{eq:sprocD2} can be factored out. Hence, the dimensions of LMIs \eqref{eq:sprocD1} and \eqref{eq:sprocD2} can be reduced by a factor of $n_\phi$, i.e., to $(Ln_\mathrm{u} + 1 + n_\mathrm{x} n_\phi)$, and to $(Ln_\mathrm{u}n_\mathrm{x} + 1 + n_\mathrm{x} n_\phi)$, respectively. The resulting LMIs are equivalent to the LMIs \eqref{eq:sprocD1} and \eqref{eq:sprocD2} before factoring. \JV{Furthermore, LMIs \eqref{eq:sprocD1}-\eqref{eq:sprocD3} are derived by invoking the matrix S-lemma \cite{vanwaarde2022noisy}. LMIs \eqref{eq:sprocD1}-\eqref{eq:sprocD3} are sufficient conditions for \eqref{eq:D1}-\eqref{eq:D3} under conditions \eqref{eq:vtilde} and \eqref{eq:tf_prop}. However, due to the imposed structural constraints due to Kronecker products, the LMIs \eqref{eq:sprocD1}-\eqref{eq:sprocD3} may, in general, be conservative.}
%	Enforcing \eqref{eq:D1}-\eqref{eq:D3} through LMIs \eqref{eq:sprocD1}-\eqref{eq:sprocD3} does not introduce any conservatism because the satisfaction of LMIs \eqref{eq:sprocD1}-\eqref{eq:sprocD3} is equivalent to the satisfaction of \eqref{eq:D1}-\eqref{eq:D3} under conditions \eqref{eq:vtilde} and \eqref{eq:tf_prop}.
\end{rem}

\JV{\begin{rem} (Computational effort) The computational effort of Problem \eqref{eq:exp_problem} scales polynomially with the state dimension $n_\mathrm{x}$, control input dimension $n_\mathrm{u}$, and the number of frequencies $L$. While this makes moderately-sized problems tractable, large-scale problems may become computationally demanding. In such cases, problem-specific structures and efficient solvers may be utilized to reduce computational demand \cite{majumdar2020recent}.\end{rem}}
The overall exploration strategy is summarized in Algorithm \ref{alg:main}. The resulting strategy optimally excites the system with inputs \eqref{eq:exploration_controller} to determine model parameters up to a user defined closeness $D_\mathrm{des}$.

\begin{algorithm}[]
\caption{Targeted exploration}
\label{alg:main}
\begin{algorithmic}[1]
\State Specify exploration length $T$, frequencies $\omega_i,\,i=1,...,L$, energy bound $\gamma_\mathrm{w}$, initial estimates $\hat{A}_0,\,\hat{B}_0$, desired accuracy of parameters $D_\mathrm{des}$.

\State Compute $\hat{V}_\phi$ and $\hat{V}_\mathrm{x}$ \eqref{eq:VhatVxhat} using the initial estimates.

\State Compute bounds $\tilde{\Gamma}_\phi$, $\tilde{\Gamma}_\mathrm{x}$ \eqref{eq:tf_prop}, and $\Gamma_\phi$, $\Gamma_\mathrm{x}$ \eqref{eq:tfy_prop} via methods described in \cite[Appendices A-C]{venkatasubramanian2023sequential}.\\ Compute matrices $\bar{W}_\phi$ \eqref{eq:PhiwPhiwH} and $\bar{W}_Z$ \eqref{eq:barWZ}.

\State Select initial candidates $\hat{Z}$ and $\hat{U}$ \eqref{eq:initcand}.

\State Set tolerance $\mathrm{tol}>0$.

\While{$\lvert \frac{\gamma_\mathrm{e} - \gamma_\mathrm{e}^*}{\gamma_\mathrm{e}}\rvert \geq \mathrm{tol}$}
\State Solve the optimization problem \eqref{eq:exp_problem}.
\State Update $\hat{Z}$ and $\hat{U}$ \eqref{eq:initcand}.
\EndWhile

\State Apply the exploration input \eqref{eq:exploration_controller} for $k=0,...,T-1$.

\State Compute parameter estimate $\hat{\theta}_T$ \eqref{eq:mean_est}; compute parameter set $\bm{\Theta}_T$ \eqref{eq:Theta_T}.

\end{algorithmic}
\end{algorithm}

\subsection{Discussion}
In what follows, we examine the key features of the proposed work and discuss its connections to the state-of-the-art. 

\textit{Summary - proposed approach:} The proposed targeted exploration strategy outlined in Algorithm \ref{alg:main} yields a multi-sine exploration input with minimal input energy to generate data from which estimates of uncertain parameters can be derived with a desired error bound \eqref{eq:exp_goal2}. The frequencies $\omega_i$ of the multi-sine input \eqref{eq:exploration_controller} are predetermined, enabling intuitive tuning based on prior knowledge about the system. The proposed approach assumes energy-bounded disturbances, as commonly considered for data-driven models~\cite{fogel1979system,vanwaarde2022noisy, berberich2022combining}. As the main result, the data-dependent uncertainty bound in Lemma \ref{lem:thetat} is utilized to derive sufficient conditions in the spectral content of the exploration inputs. The proposed exploration is \emph{targeted}, as the optimized amplitudes at different frequencies impact both the magnitude and the shape/orientation of the remaining uncertainty after exploration. Furthermore, the proposed exploration strategy is robust, i.e., parametric uncertainty is accounted for by using robust control tools.

\textit{Limitations:} The proposed exploration strategy relies on the energy-bound of the disturbances $\gamma_\mathrm{w}$. The strategy requires solving an SDP iteratively to mitigate suboptimality arising from the convex relaxation procedure. Additionally, the proposed robust exploration strategy is more conservative for large initial uncertainty.

\textit{Related works:} The derived targeted exploration strategy is similar to, and inspired by \cite{venkatasubramanian2023sequential, umenberger2019robust, ferizbegovic2019learning, barenthin2008identification}, however, with a few crucial differences. A targeted exploration method is proposed in \cite{barenthin2008identification} to identify parameters up to a desired accuracy. However, the conditions are not robust to uncertainty, requiring iterative experiments in practice. In \cite{ferizbegovic2019learning}, the proposed targeted exploration strategy assumes independent data, and hence lacks applicability to correlated time-series data from a single trajectory. Furthermore, the methods in \cite{ferizbegovic2019learning} and \cite{umenberger2019robust} do not yield any guarantees for exploration since the uncertainty bounds are heuristically approximated. \JV{The method in \cite{venkatasubramanian2023sequential} robustly accounts for parametric uncertainties and provides an \textit{a priori} guaranteed bound on the uncertainty after exploration. However, all the methods \cite{barenthin2008identification,ferizbegovic2019learning,umenberger2019robust,venkatasubramanian2023sequential} assume i.i.d. Gaussian disturbances. In contrast to these methods, we consider energy-bounded disturbances without assumptions on the distribution or independence of the disturbances. This encompasses a broader class of uncertainties, including those arising from unmodeled dynamics or nonlinearities, and enables the development of a targeted exploration strategy with guarantees. In our proposed targeted exploration method, we quantify and guarantee \textit{a priori} uncertainty bounds on the parameters obtained from a single experiment by robustly accounting for the impact of the uncertain model parameters and disturbances, similar to \cite{venkatasubramanian2023sequential}. A notable difference to the theoretical analysis in \cite{venkatasubramanian2023sequential} is the appearance of the term $Z Z^\top$ in Theorem \ref{thm:sufficientcond1}, which cannot be directly addressed with standard robust control tools and presented significant additional challenges in this work.}

\textit{Application to dual control:} The proposed targeted exploration strategy can be integrated into a robust \textit{dual control} framework to guarantee a desired closed-loop performance after exploration \cite{barenthin2008identification, venkatasubramanian2023sequential}. This can be achieved by co-designing the targeted exploration problem \eqref{eq:exp_problem} with a robust gain-scheduled controller that accounts for changes in uncertainty during exploration \cite{venkatasubramanian2023sequential}. Such a joint design allows optimizing closed-loop performance after exploration, while ensuring exploration reduces only the necessary uncertainty to achieve the desired performance within exploration energy constraints.

%\section{Numerical Example}\label{sec:numerical}
%\input{sections/numericalsimulations}

\section{Numerical Example}\label{sec:numerical}
In this section, we demonstrate the applicability of the proposed targeted exploration strategy to a nonlinear system using a numerical example. Numerical simulations\footnote{\scriptsize{The source code for the simulations is available at https:$//$github.com$/$jananivenkatasubramanian$/$NonstochTE}} were performed on MATLAB using CVX \cite{cvx} in conjunction with the solver MOSEK.

\textit{Problem Setup:}
We consider a chain of two mass-spring-damper systems. The model equations are given by
\vspace{-20pt}
\begin{small}
\begin{align}
\nonumber
    m_1 \ddot{p}_1 &= -(k_1 + k_2)p_1 - (d_1+d_2)\dot{p}_1 + k_2 p_2 + d_2 \dot{p}_2 + F_1^\mathrm{nl}, \\
\label{eq:modeleq}    
    m_2 \ddot{p}_2 &= k_2 p_1 + d_2 \dot{p}_1 - k_2 p_2 - d_2 \dot{p}_2 + F_2^\mathrm{nl}+u,\\
\label{eq:coulombfriction}    
    F_i^\mathrm{nl}&= \zeta_i \tanh (\alpha_i \dot{p}_i),
\end{align}
\end{small}
with positions $p_i$, velocities $\dot{p}_i$, masses $m_i > 0$, spring constants $k_i \geq 0$, damping coefficients $d_i \geq 0$, nonlinear Coulomb frictions $F_i^\mathrm{nl}$ with constants $\alpha_i,\, \zeta_i \geq 0$, $i \in \{1,2\}$, and input $u \in \mathbb{R}$. We consider $x=[p_1,\dot{p}_1, p_2, \dot{p}_2]^\top \in \mathbb{R}^4$ and an Euler discretization with sampling period $T_s=0.5$. The resulting system dynamics correspond to our setup \eqref{eq:sys} with the disturbances $w_k=\left[0,\frac{F_1^\mathrm{nl}(x_2)}{m_1},0,\frac{F_1^\mathrm{nl}(x_4)}{m_2}\right]^\top$. The disturbances satisfy Assumption \ref{a0} with $\gamma_\mathrm{w}=T\left(\frac{\zeta_1^2}{m_1^2}+\frac{\zeta_2^2}{m_2^2}\right)$. In this example, we consider the following true values for the masses, spring constants and damping coefficients: $m_1=1$, $m_2=2$, $k_1=1$, $k_2=1.5$, $d_1=0.5$ and $d_2=1.1$ \eqref{eq:modeleq}. We select the coefficients of Coulomb friction as $\alpha_1=1$ and $\alpha_2=1$ in \eqref{eq:coulombfriction}. Later, we vary $\zeta_1$ and $\zeta_2$ to study the effect of the disturbance bound $\gamma_\mathrm{w}$.

The goal of the proposed targeted exploration strategy is to achieve a desired error bound $D_\mathrm{des}^{-1}=10^{0} I_{n_\phi}$ on the parameters \eqref{eq:exp_goal2}. Furthermore, we select $L=20$ equally-spaced frequencies $\omega_i\in\{0,0.05,0.1,0.15,...,0.95\}$. We set $T=100$ and $\epsilon=0.5$. Throughout the example, the constants $\Gamma_\phi$, $\Gamma_\mathrm{x}$, $\tilde{\Gamma}_\phi$ and $\tilde{\Gamma}_\mathrm{x}$ are computed using the scenario approach \cite[Appendix C]{venkatasubramanian2023sequential} with confidence level $\beta=10^{-10}$. In what follows, we analyse the effectiveness and conservatism of the proposed targeted exploration strategy.

\textit{Required input energy $\gamma_\mathrm{e}^2$ for different energy bounds $\gamma_\mathrm{w}$:} To study the effect of the energy bound $\gamma_\mathrm{w}$ on the required input energy $\gamma_\mathrm{e}^2$, we select the initial uncertainty level as $D_0^{-1}=10^{-4}I_{n_\phi}$ and the initial estimate as $\hat{\theta}_0=\theta_\mathrm{tr}+\frac{\|D_0\|^{\frac{1}{2}}}{\|\theta_\mathrm{tr}\|}\theta_\mathrm{tr}$, i.e., on the boundary of $\mathbf{\Theta}_0$ (cf. Assumption \ref{a1}). We run five trials for the following energy bounds $\gamma_\mathrm{w} \in \{10^{-4},10^{-3},10^{-2},10^{-1},10^{0}\}$. Each trial comprises: $\mathrm{(i)}$ computing the corresponding constants $\Gamma_\phi$, $\Gamma_\mathrm{x}$, $\tilde{\Gamma}_\phi$ and $\tilde{\Gamma}_\mathrm{x}$, $\mathrm{(ii)}$ executing Algorithm \ref{alg:main} to obtain the exploration inputs \eqref{eq:exploration_controller} and the required input energy $\gamma_\mathrm{e}^2$.
\begin{figure}[t!]
\begin{center}
\includegraphics[width=0.5\textwidth]{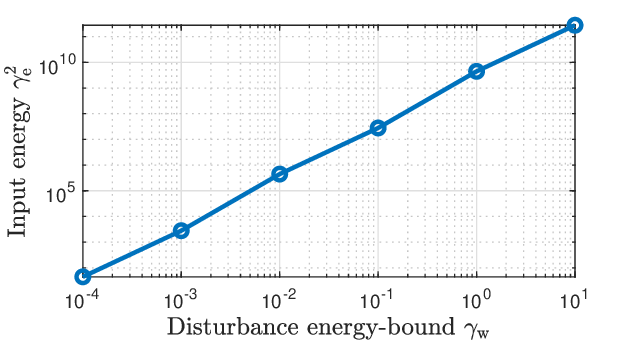}
\end{center}
\caption{Illustration of the exploration input energy $\gamma_\mathrm{e}^2$, in comparison with the disturbance energy bound $\gamma_\mathrm{w}$ for the initial uncertainty level $\|D_0\|=10^2$.}
\label{fig:ge-gw} 
\end{figure}
From Figure \ref{fig:ge-gw}, it can be observed that the input energy $\gamma_\mathrm{e}^2$ scales roughly linearly with the disturbance energy bound $\gamma_\mathrm{w}$. In general, $\gamma_\mathrm{e} \to 0$ as $\gamma_\mathrm{w} \to 0$, i.e., the exploration input energy $\gamma_\mathrm{e}$ reduces to zero as the energy of the disturbances $\gamma_\mathrm{w}$ reduce to zero.
%Furthermore, the input energy needs to be significantly larger than the disturbance energy to achieve the desired accuracy of the estimated parameters. 

\textit{\JV{Conservatism related to the initial uncertainty bound $D_0$ and comparison with naive exploration:}} To study the effect of the initial uncertainty bound $D_0^{-1}$, we select eight uncertainty levels $D_0 \in \{10^i I_{n_\phi}, 5 \times 10^i I_{n_\phi}\}$, $i = 4,...,7$ and $\gamma_\mathrm{w}=1$. Based on $\gamma_\mathrm{w}$, we compute $\zeta_1 = \zeta_2= \frac{m_1 m_2 \sqrt{\gamma_\mathrm{w}}}{\sqrt{T(m_1^2 + m_2^2)}}$. \JV{Furthermore, to highlight the benefits of the proposed method, we compare the targeted exploration method with a naive exploration strategy that uses inputs of equal total energy but with non-optimized amplitudes, i.e., inputs where the energy is uniformly distributed across all frequencies. This baseline does not rely on any model information and illustrates how optimization improves input design.} Each trial comprises: $\mathrm{(i)}$ computing the initial estimate as $\hat{\theta}_0=\theta_\mathrm{tr}+\frac{\|D_0\|^{\frac{1}{2}}}{\|\theta_\mathrm{tr}\|}\theta_\mathrm{tr}$, $\mathrm{(ii)}$ computing the corresponding constants $\Gamma_\phi$, $\Gamma_\mathrm{x}$, $\tilde{\Gamma}_\phi$ and $\tilde{\Gamma}_\mathrm{x}$, $\mathrm{(iii)}$ executing Algorithm~\ref{alg:main} to obtain the exploration inputs \eqref{eq:exploration_controller}, and generate a dataset by applying the targeted exploration input, \JV{$\mathrm{(iv)}$ computing the naive exploration inputs, and generating a corresponding dataset.} From both datasets, we compute the corresponding estimate of the parameters $\hat{\theta}_T$, and $\|G \cdot P\|$, the error bound guaranteed by Lemma \ref{lem:thetat} with the data-dependent covariance matrix $P$ \eqref{eq:covar_est} and scaling $G$ \eqref{eq:G}.
\begin{figure}[t!]
\begin{center}
\includegraphics[width=0.5\textwidth]{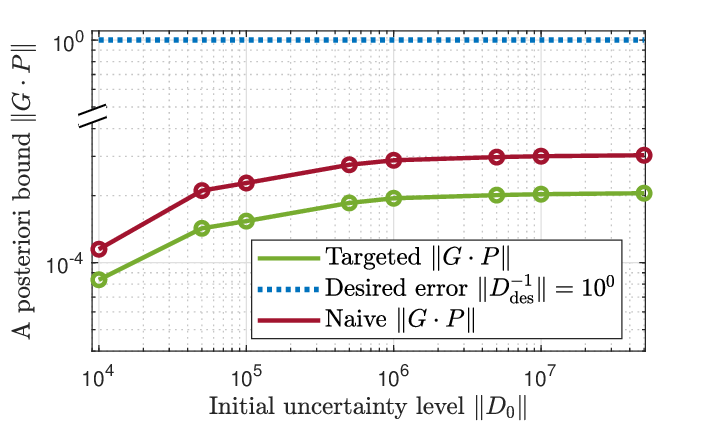}
\end{center}
\caption{Illustration of the \textit{a posteriori} guaranteed bound on the squared error of the parameters $\|G \cdot P\|$ for both targeted and naive exploration with same input energy, in comparison with the desired bound on the squared error $\|D_\mathrm{des}^{-1}\|$ for different initial uncertainty bounds $D_0^{-1}$.}
\label{fig:gw-D0-error} 
\end{figure}
Figure \ref{fig:gw-D0-error} shows that the targeted exploration inputs achieve the exploration goal for all tested initial uncertainty levels $D_0$ and energy bound $\gamma_\mathrm{w}$, i.e., the guaranteed error bound on the parameters $\|G \cdot P \|$ is lower than the desired error bound $\|D_\mathrm{des}^{-1}\|~=~10^{0}$. The guaranteed bound is significantly lower than the desired bound due to the inherent conservatism of the proposed strategy which is derived by utilizing worst case bounds on the transfer matrices $\Gamma_\phi$, $\Gamma_\mathrm{x}$, $\tilde{\Gamma}_\phi$ and $\tilde{\Gamma}_\mathrm{x}$, and worst case bounds $\Bar{W}_\phi$, $\bar{W}_Z$. This conservatism decreases as the initial uncertainty is reduced, or equivalently, as $D_0$ increases, since smaller initial uncertainty results in tighter and more accurate bounds on $\Gamma_\phi$, $\Gamma_\mathrm{x}$, $\tilde{\Gamma}_\phi$ and $\tilde{\Gamma}_\mathrm{x}$, yielding more efficient exploration. \JV{Furthermore, Figure \ref{fig:gw-D0-error} shows that the proposed targeted exploration method guarantees a lower error bound $\|G\cdot P\|$, by roughly $50\%$, compared to naive exploration under the same energy budget, thereby highlighting the benefits of our approach.}

\textit{Sensitivity related to the initial estimate $\hat{\theta}_0$:} To study the sensitivity related to the initial estimate $\hat{\theta}_0$, we select three different uncertainty levels $D_0 \in \{ 10^i I_{n_\phi} \}$, $i = 4,5,6$, with $\gamma_\mathrm{w}=1$. For each level, we run 50 trials. Each trial comprises: $\mathrm{(i)}$ generating a random initial estimate $\hat{\theta}_0$ that satisfies Assumption \ref{a1}, i.e., $(\hat{\theta}_0-\theta_\mathrm{tr})^\top (D_0 \otimes I_{n_\mathrm{x}})  (\hat{\theta}_0-\theta_\mathrm{tr}) \leq 1$, $\mathrm{(ii)}$ computing constants $\Gamma_\phi$, $\Gamma_\mathrm{x}$, $\tilde{\Gamma}_\phi$ and $\tilde{\Gamma}_\mathrm{x}$, $\mathrm{(iii)}$ executing Algorithm \ref{alg:main} to obtain the exploration inputs \eqref{eq:exploration_controller} and the input energy $\gamma_\mathrm{e}^2$.
\begin{figure}[t!]
	\begin{center}
		\includegraphics[width=0.49\textwidth]{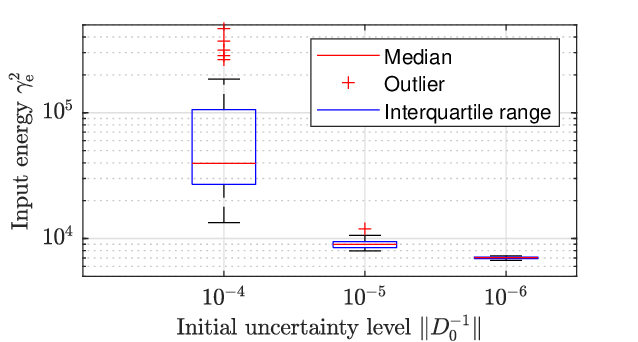}
	\end{center}
	\caption{Illustration of the distribution of the input energy $\gamma_\mathrm{e}^2$ for different initial uncertainty bounds $D_0^{-1}$.}
	\label{fig:sensitivity1} 
\end{figure}
Figure \ref{fig:sensitivity1} shows that the variability of input energy $\gamma_\mathrm{e}^2$ is higher for larger initial uncertainty levels $\|D_0^{-1}\|$.  This variability significantly reduces as initial uncertainty levels reduce. This behaviour arises because smaller initial uncertainty levels yield tighter bounds on $\Gamma_\phi$, $\Gamma_\mathrm{x}$, $\tilde{\Gamma}_\phi$ and $\tilde{\Gamma}_\mathrm{x}$, resulting in more consistent input energy.

\textit{Sensitivity related to the desired error bound $\|D_\mathrm{des}\|$:} To study the sensitivity related to the desired uncertainty level $\|D_\mathrm{des}\|$, we select six levels $D_\mathrm{des}=10^{i} I_{n_\phi}$, $i \in \{-4,-3,-2,-1,0,1\}$. We set the initial uncertainty level to $D_0=10^4 I_{n_\phi}$, $\gamma_\mathrm{w}=1$, choose $\hat{\theta}_0=\theta_\mathrm{tr}+\frac{\|D_0\|^{\frac{1}{2}}}{\|\theta_\mathrm{tr}\|}\theta_\mathrm{tr}$, and compute the corresponding constants $\Gamma_\phi$, $\Gamma_\mathrm{x}$, $\tilde{\Gamma}_\phi$ and $\tilde{\Gamma}_\mathrm{x}$. For each $\|D_\mathrm{des}\|$, we execute Algorithm \ref{alg:main} to obtain the exploration inputs \eqref{eq:exploration_controller} and the input energy $\gamma_\mathrm{e}^2$.
\begin{figure}[t!]
	\begin{center}
		\includegraphics[width=0.49\textwidth]{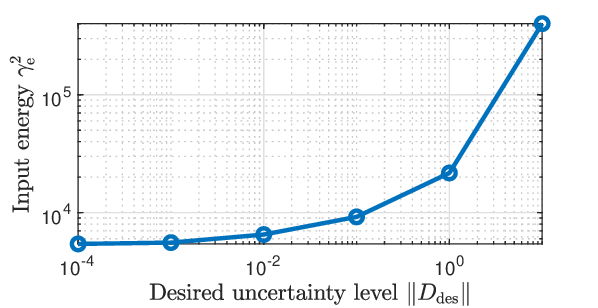}
	\end{center}
	\caption{Illustration of input energy $\gamma_\mathrm{e}^2$ for different desired uncertainty levels $\|D_\mathrm{des}\|$.}
	\label{fig:ddes} 
\end{figure}
\JV{Figure \ref{fig:ddes} shows that achieving a small uncertainty (large $\|D_\mathrm{des}\|$) requires higher input energy. In contrast, the required input energy scales linearly with $\|D_\mathrm{des}\|$ in the stochastic disturbances setting \cite{venkatasubramanian2023sequential}. This difference is due to the data-dependent scaling $G$ \eqref{eq:G} of the non-falsified set $\mathbf{\Theta}_T$ (cf. Lemma \ref{lem:thetat}).}

%Furthermore, in order to shrink the ellipsoid for smaller desired uncertainty levels, shrinking the scaling $G$ \eqref{eq:G} requires increasingly large $X^\top X$ relative to $\gamma_\mathrm{w}$, and thereby increasingly large input energy $\gamma_\mathrm{e}^2$. This leads to a superlinear relationship between the required input energy $\gamma_\mathrm{e}^2$ and the desired uncertainty level $\|D_\mathrm{des}\|$. Consequently, designing exploration inputs for large values of $D_\mathrm{des}$, i.e., small desired uncertainty in the parameters, quickly becomes numerically infeasible since the required input energy may grow prohibitively large.

\ifbool{arxivversion}{\textit{Computational time:} The average execution time for solving Problem \eqref{eq:exp_problem} over 10 trials is approximately 45 seconds. The simulations were carried out on a system with an AMD Ryzen 7 5700U processor and 16.0 GB RAM.}{}
%In summary, the simulation results demonstrate the applicability of the proposed exploration strategy to systems with unmodeled nonlinearities. Given an initial estimate $\hat{\theta}_0$ and initial uncertainty bound $D_0^{-1}$, the targeted exploration strategy guarantees an \textit{a priori} error bound on parameters estimated after exploration. 

\section{Conclusion}\label{sec:conclusion}
In this article, we present a targeted exploration strategy for linear systems subject to energy-bounded disturbances. We derived LMIs that can robustly guarantee an $\textit{a priori}$ error-bound on the estimated parameters after exploration. The proposed strategy utilizes multi-sine inputs in selected frequencies and optimized amplitudes to shape the uncertainty bound over the parameters in a targeted manner. We demonstrate the applicability of the proposed targeted exploration strategy to systems with nonlinearities, while analysing its benefits and inherent conservatism, with a numerical example. This is the first targeted exploration approach that robustly ensures a user-chosen accuracy on the parameters without requiring any independence conditions on the disturbances. \JV{We view addressing the limitations related to conservatism and scalability as an important direction for future work.}

% We derive sufficient conditions on the spectral content of the exploration data  Using these sufficient conditions, we derive robust LMIs for exploration while accounting for initial parametric uncertainty.  
% In summary, to the best of the authors' knowledge, the proposed approach is the first principled, computationally tractable, and robust finite-time targeted exploration (optimal experiment design) strategy for linear systems with energy-bounded disturbances without required independence of disturbances.

\bibliographystyle{plain}       
\bibliography{lit}          

\appendix

\ifbool{arxivversion}{\section{Proof of Lemma \ref{lemma:ub_Zlb_tight}}\label{app:spectralbounds}
	\begin{pf} 
\begin{table*}
\begin{align}\label{eq:ZZtop_ineq}
\nonumber
    Z Z^\top &
    \overset{\eqref{eq:spectral_amplitude}}{=} T \begin{bmatrix}
        D_\mathrm{des}^{\frac{1}{2} \top} \left( \left(\sum_{i=1}^T \bar{x}(\omega_i)^\mathsf{H} \bar{x}(\omega_i)\right) \otimes I_{n_\phi} \right)D_\mathrm{des}^{\frac{1}{2}} & D_\mathrm{des}^{\frac{1}{2} \top} \left( \sum_{i=1}^{T}\bar{x}(\omega_i)^\top \left(\bar{\phi}(\omega_i)^\mathsf{H} \otimes I_{n_\mathrm{x}}\right) \right)\otimes I_{n_\phi}\\
        \left( \sum_{i=1}^{T}\bar{x}(\omega_i)^\top \left(\bar{\phi}(\omega_i)^\mathsf{H} \otimes I_{n_\mathrm{x}}\right) \right)^\mathsf{H}\otimes I_{n_\phi}  D_\mathrm{des}^{\frac{1}{2}} & \left( \left( \sum_{i=1}^{T} \bar{\phi}(\omega_i) \bar{\phi}(\omega_i)^\mathsf{H} \right) \otimes I_{n_\mathrm{x}}\right)\otimes I_{n_\phi}
    \end{bmatrix}\\
    & \succeq T \begin{bmatrix}
        D_\mathrm{des}^{\frac{1}{2} \top} \left( \left(\sum_{i=1}^L \bar{x}(\omega_i)^\mathsf{H} \bar{x}(\omega_i)\right) \otimes I_{n_\phi} \right)D_\mathrm{des}^{\frac{1}{2}} & D_\mathrm{des}^{\frac{1}{2} \top} \left( \sum_{i=1}^{L}\bar{x}(\omega_i)^\top \left(\bar{\phi}(\omega_i)^\mathsf{H} \otimes I_{n_\mathrm{x}}\right) \right)\otimes I_{n_\phi}\\
        \left( \sum_{i=1}^{L}\bar{x}(\omega_i)^\top \left(\bar{\phi}(\omega_i)^\mathsf{H} \otimes I_{n_\mathrm{x}}\right) \right)^\mathsf{H}\otimes I_{n_\phi}  D_\mathrm{des}^{\frac{1}{2}} & \left( \left( \sum_{i=1}^{L} \bar{\phi}(\omega_i) \bar{\phi}(\omega_i)^\mathsf{H} \right) \otimes I_{n_\mathrm{x}}\right)\otimes I_{n_\phi}
    \end{bmatrix}
\end{align}
\end{table*}
The proof is divided into two parts. We first derive \eqref{eq:phiphitight} and then derive \eqref{eq:ZZtight}.

\textbf{Part I.} From the Parseval-Plancherel identity, we have
\begin{align}\label{eq:parsphiphi}
\nonumber
    \Phi \Phi^\top=\sum_{k=0}^{T-1} \phi_k \phi_k^\top &\overset{\eqref{eq:spectral_amplitude}}{=}T \left( \sum_{i=1}^{T} \bar{\phi}(\omega_i) \bar{\phi}(\omega_i)^\mathsf{H} \right)\\
\nonumber
    &\succeq T \left( \sum_{i=1}^{L} \bar{\phi}(\omega_i) \bar{\phi}(\omega_i)^\mathsf{H} \right)\\
    &=T \left( \bar{\Phi}_\mathrm{u}+\bar{\Phi}_\mathrm{w} \right) \left( \bar{\Phi}_\mathrm{u}+\bar{\Phi}_\mathrm{w} \right)^\mathsf{H}.
\end{align}
By Young's inequality \cite{caverly2019lmi}, for any $\epsilon >0$, we have
\begin{align}
\bar{\Phi}_\mathrm{u} \bar{\Phi}_\mathrm{w}^\mathsf{H} + \bar{\Phi}_\mathrm{w} \bar{\Phi}_\mathrm{u}^\mathsf{H} \succeq -\epsilon \bar{\Phi}_\mathrm{u} \bar{\Phi}_\mathrm{u}^\mathsf{H} - \frac{1}{\epsilon} \bar{\Phi}_\mathrm{w} \bar{\Phi}_\mathrm{w}^\mathsf{H}
\end{align}
and hence,
\begin{equation}\label{eq:yi}
(\bar{\Phi}_\mathrm{u}+\bar{\Phi}_\mathrm{w})(\bar{\Phi}_\mathrm{u}+\bar{\Phi}_\mathrm{w})^\mathsf{H} \succeq (1-\epsilon) \bar{\Phi}_\mathrm{u} \bar{\Phi}_\mathrm{u}^\mathsf{H} - \left(\frac{1-\epsilon}{\epsilon}\right) \bar{\Phi}_\mathrm{w} \bar{\Phi}_\mathrm{w}^\mathsf{H}.
\end{equation}
By inserting Inequality \eqref{eq:yi} in Inequality \eqref{eq:parsphiphi}, we get \eqref{eq:phiphitight}. \\
\textbf{Part II.} 
From the Parseval-Plancheral identity, we have \eqref{eq:ZZtop_ineq}, which can be written as
\begin{align}\label{eq:ZZtop}
    Z Z^\top \succeq T \bar{Z} \bar{Z}^\mathsf{H} = T \left(\bar{Z}_\mathrm{u}+\bar{Z}_\mathrm{w} \right)\left(\bar{Z}_\mathrm{u}+\bar{Z}_\mathrm{w} \right)^\mathsf{H}.
\end{align}
By Young's inequality \cite{caverly2019lmi}, for any $\epsilon >0$, we have
\begin{align}
\bar{Z}_\mathrm{u} \bar{Z}_\mathrm{w}^\mathsf{H} + \bar{Z}_\mathrm{w} \bar{Z}_\mathrm{u}^\mathsf{H} \succeq -\epsilon \bar{Z}_\mathrm{u} \bar{Z}_\mathrm{u}^\mathsf{H} - \frac{1}{\epsilon} \bar{Z}_\mathrm{w} \bar{Z}_\mathrm{w}^\mathsf{H}
\end{align}
and hence,
\begin{equation}\label{eq:yi2}
(\bar{Z}_\mathrm{u}+\bar{Z}_\mathrm{w})(\bar{Z}_\mathrm{u}+\bar{Z}_\mathrm{w})^\mathsf{H} \succeq (1-\epsilon) \bar{Z}_\mathrm{u} \bar{Z}_\mathrm{u}^\mathsf{H} - \left(\frac{1-\epsilon}{\epsilon}\right) \bar{Z}_\mathrm{w} \bar{Z}_\mathrm{w}^\mathsf{H}.
\end{equation}
By inserting Inequality \eqref{eq:yi2} in Inequality \eqref{eq:ZZtop}, we get \eqref{eq:ZZtight}. $\hfill\square$
\end{pf}

}{}

\ifbool{arxivversion}{\section{Proof of Lemma \ref{lem:noisebounds}}\label{app:noisebounds}
\begin{pf} The proof is provided in two parts. In the first part, we prove \eqref{eq:PhiwPhiwH}, and in the second part, we prove \eqref{eq:barWZ}.

\textbf{Part I.} In order to determine a bound on $\bar{\Phi}_\mathrm{w}$ \eqref{eq:phibarcompact} of the form in \eqref{eq:PhiwPhiwH}, we first determine a bound on $W$ \eqref{eq:VxYxW}. In particular, we have
\begin{align}\label{eq:Wnorm}
\|W\| = \max_{i=1,\dots,L}\| \bar{w}(\omega_i)\|
\end{align}
and
\begin{align}\label{eq:wi_sum}
    \max_{i=1,\dots,L}\| \bar{w}(\omega_i) \|^2 \leq \sum_{i=1}^L \| \bar{w}(\omega_i) \|^2 \leq \sum_{\omega_i \in \Omega_T} \| \bar{w}(\omega_i) \|^2.
\end{align}
From the Parseval-Plancherel identity and the energy-bound on the disturbance \eqref{eq:noisebound}, we have
\begin{align}\label{eq:wmaxnorm}
    \sum_{\omega_i \in \Omega_T} \| \bar{w}(\omega_i) \|^2 = \frac{1}{T} \sum_{i=0}^{T-1}\|w_k\|^2 \overset{\eqref{eq:noisebound}}{\leq} \frac{\gamma_\mathrm{w}}{T},
\end{align}
and hence, from \eqref{eq:Wnorm} and \eqref{eq:wmaxnorm}, we have
\begin{align}\label{eq:normw}
    \|W\| \leq \sqrt{\frac{\gamma_\mathrm{w}}{T}} \text{ and } W W^\mathsf{H} \preceq \frac{\gamma_\mathrm{w}}{T} I_{n_\mathrm{x}L}.
\end{align}
Starting from \eqref{eq:phibarcompact}, we get
\begin{align*}
\nonumber
    \bar{\Phi}_\mathrm{w} \bar{\Phi}_\mathrm{w}^\mathsf{H} & = Y_{\phi,\mathrm{tr}} W W^\mathsf{H} Y_{\phi,\mathrm{tr}}^\mathsf{H}\\
\nonumber
    & \overset{\eqref{eq:normw}}{\preceq}  Y_{\phi,\mathrm{tr}} \left( \frac{\gamma_\mathrm{w}}{T} I_{n_\mathrm{x}L}  \right) Y_{\phi,\mathrm{tr}}^\mathsf{H}\\
    & \overset{\eqref{eq:tfy_prop}}{\preceq} \frac{\gamma_\mathrm{w}}{T} \Gamma_\phi = \bar{W}_\phi. \tag{cf. \eqref{eq:PhiwPhiwH}}
\end{align*}
\textbf{Part II.} From \eqref{eq:PhiwPhiwH}, we have
\begin{align}\label{eq:Phiwnorm}
    \|\bar{\Phi}_\mathrm{w}\| \leq  \sqrt{\frac{\gamma_\mathrm{w}}{T}} \|\Gamma_\phi\|^{\frac{1}{2}}.
\end{align}
Similarly, from \eqref{eq:xbarcompact}, we have
\begin{align}\label{eq:Xbarnorm}
    \|\bar{X}_\mathrm{w}\| \leq & \|Y_\mathrm{\mathrm{x},tr}\| \|W\mathbf{1}_L\|\\
    \leq & \|Y_\mathrm{\mathrm{x},tr}\|  \max_{i=1,\dots,L}\| \bar{w}(\omega_i) \|^2 \overset{\eqref{eq:tfy_prop}, \atop \eqref{eq:wi_sum}}{=}
    \sqrt{\frac{\gamma_\mathrm{w}}{T}} \|\Gamma_\mathrm{x}\|^{\frac{1}{2}}.\nonumber
\end{align}
Recall that $\bar{Z}_\mathrm{w}=\begin{bmatrix}
        D_\mathrm{des}^{\frac{1}{2}\top} \left(\bar{X}_\mathrm{w}^\mathsf{H}  \otimes I_{n_\phi}\right)\\
        \bar{\Phi}_\mathrm{w} \otimes I_{n_\mathrm{x}n_\phi}
    \end{bmatrix}$ from \eqref{eq:Zspectral}.
    Note that
\begin{align}\label{eq:wzHwz}
\nonumber
   \bar{Z}_\mathrm{w}^\mathsf{H} \bar{Z}_\mathrm{w} =&  \left(\bar{X}_\mathrm{w}  \otimes I_{n_\phi}\right)D_\mathrm{des}^{\frac{1}{2}}D_\mathrm{des}^{\frac{1}{2}\top}\left(\bar{X}_\mathrm{w}^\mathsf{H}  \otimes I_{n_\phi}\right)\\
   &+ \left(\bar{\Phi}_\mathrm{w}^\mathsf{H} \bar{\Phi}_\mathrm{w}\otimes I_{n_\mathrm{x} n_\phi}\right).
\end{align}
By taking the induced norm of both sides of \eqref{eq:wzHwz}, followed by the application of the triangle inequality, we have
\begin{align}\label{eq:Zbarnorm}
\nonumber
    \| \bar{Z}_\mathrm{w} \|^2 \leq & \|D_\mathrm{des}^{\frac{1}{2}\top}(\bar{X}_\mathrm{w}^\mathsf{H}  \otimes I_{n_\mathrm{x}n_\phi})\|^2+\|\bar{\Phi}_\mathrm{w}\|^2\\
    \overset{\eqref{eq:Phiwnorm},\atop \eqref{eq:Xbarnorm}}{\leq} & \left(\frac{\gamma_\mathrm{w}}{T}\|\Gamma_\mathrm{x}\|\|D_\mathrm{des}\| +     \frac{\gamma_\mathrm{w}}{T}\|\Gamma_\phi\| \right).
\end{align}
Finally, since $\|\bar{Z}_\mathrm{w} \bar{Z}_\mathrm{w}^\mathsf{H}\|=\|\bar{Z}_\mathrm{w}^\mathsf{H} \bar{Z}_\mathrm{w}\|=\|\bar{Z}_\mathrm{w}\|^2$, from \eqref{eq:Zbarnorm} we have
 \begin{align*}
 \nonumber
    \bar{Z}_\mathrm{w} \bar{Z}_\mathrm{w}^\mathsf{H} & \preceq \|\bar{Z}_\mathrm{w}\|^2 I_{(n_\phi + n_\mathrm{x} n_\phi^2)}\\
    \nonumber
    &= \left(\frac{\gamma_\mathrm{w}}{T}\|\Gamma_\mathrm{x}\|\|D_\mathrm{des}\| +     \frac{\gamma_\mathrm{w}}{T}\|\Gamma_\phi\| \right)I_{(n_\phi + n_\mathrm{x} n_\phi^2)}
\end{align*}
which yields \eqref{eq:barWZ}. $\hfill\square$
\end{pf}}{}

\ifbool{arxivversion}{\section{Proof of Proposition \ref{prop2}}\label{app:prop2}  
	\begin{pf}
	Starting from Inequality \eqref{eq:prop2}, we have
	\begin{align*}
		\nonumber
		0 \preceq & \begin{bmatrix}
			T \left( (1-\epsilon) \bar{\Phi}_\mathrm{u} \bar{\Phi}_\mathrm{u}^\mathsf{H} - \left(\tfrac{1-\epsilon}{\epsilon}\right) \bar{W}_\phi \right)  - \gamma_\mathrm{w} D_\mathrm{des} & 0\\0 & 0
		\end{bmatrix}\\
		\nonumber
		& + T\left((1-\epsilon) \bar{Z}_\mathrm{u} \bar{Z}_\mathrm{u}^\mathsf{H} - \left( \tfrac{1-\epsilon}{\epsilon}\right) \bar{W}_Z\right)\\
		\nonumber
		\overset{\eqref{eq:PhiwPhiwH},\atop \eqref{eq:barWZ}}{\preceq} & \begin{bmatrix}
			T\left( (1-\epsilon) \bar{\Phi}_\mathrm{u} \bar{\Phi}_\mathrm{u}^\mathsf{H} - \left(\tfrac{1-\epsilon}{\epsilon}\right) \bar{\Phi}_\mathrm{w} \bar{\Phi}_\mathrm{w}^\mathsf{H}\right) - \gamma_\mathrm{w} D_\mathrm{des} & 0\\0 & 0
		\end{bmatrix}\\
		\nonumber
		& + T\left( (1-\epsilon) \bar{Z}_\mathrm{u} \bar{Z}_\mathrm{u}^\mathsf{H} - \left( \tfrac{1-\epsilon}{\epsilon}\right) \bar{Z}_\mathrm{w} \bar{Z}_\mathrm{w}^\mathsf{H}\right)\\
		\overset{\eqref{eq:phiphitight}, \atop \eqref{eq:ZZtight}}{\preceq} &\begin{bmatrix}
			(\Phi \Phi^\top)\otimes I_{n_\mathrm{x}} - \gamma_\mathrm{w} D_\mathrm{des} & 0\\  
			0 & 0
		\end{bmatrix} + Z Z^\top.
	\end{align*}
	The condition \eqref{eq:prop2} corresponds to the condition \eqref{eq:pe5} in Theorem \ref{thm:sufficientcond1}, and hence the exploration goal \eqref{eq:exp_goal2} is achieved. $\hfill\square$
\end{pf}}{}

\ifbool{arxivversion}{\section{Proof of Proposition \ref{prop4}}\label{app:convexrel}\begin{pf}
	For matrices $\bar{Z}_\mathrm{u}$ and $\hat{Z}$, we have $\bar{Z}_\mathrm{u}\bar{Z}_\mathrm{u}^\mathsf{H} - \bar{Z}_\mathrm{u}\hat{Z}^\mathsf{H}- \hat{Z} \bar{Z}_\mathrm{u}^\mathsf{H} + \hat{Z} \hat{Z}^\mathsf{H} = (\bar{Z}_\mathrm{u} - \hat{Z})(\bar{Z}_\mathrm{u} - \hat{Z})^\mathsf{H}\succeq 0$, and hence, we have
	\begin{align}\label{eq:zzconvex}
		\bar{Z}_\mathrm{u}\bar{Z}_\mathrm{u}^\mathsf{H} \succeq \bar{Z}_\mathrm{u}\hat{Z}^\mathsf{H}+ \hat{Z} \bar{Z}_\mathrm{u}^\mathsf{H} -\hat{Z} \hat{Z}^\mathsf{H}.
	\end{align}
	%From Lemma~\ref{lem:convexrel1}, we have
	
	Inserting Inequality \eqref{eq:zzconvex} in Inequality \eqref{eq:prop3}, and multiplying the resulting inequality by $T$ yields Inequality \eqref{eq:prop2}. Hence, if there exists matrices $\bar{\Phi}_\mathrm{u}$, $\bar{Z}_\mathrm{u}$ and $\hat{Z}$ that satisfy \eqref{eq:prop3}, then the condition in Proposition \ref{prop2} is satisfied and the exploration goal \eqref{eq:exp_goal2} is achieved. $\hfill\square$
\end{pf}}{}

\ifbool{arxivversion}{\section{Proof of Theorem \ref{thm:finalthm}}\label{app:finalthm}  
\begin{pf}
Inequalities \eqref{eq:D1}-\eqref{eq:D3} imply the exploration goal \eqref{eq:exp_goal2} due to Proposition \ref{prop4}. In what follows, we prove that Inequalities \eqref{eq:sprocD1}-\eqref{eq:sprocD3} imply Inequalities \eqref{eq:D1}-\eqref{eq:D3}, respectively. In particular, we utilize the matrix S-lemma \cite{boyd2004convex,vanwaarde2022noisy} to account for uncertainties in $V_{\mathrm{x,tr}}$, $V_{\mathrm{\phi,tr}}$, $Y_{\mathrm{x,tr}}$, and $Y_{\mathrm{\phi,tr}}$ satisfying bounds \eqref{eq:tf_prop}, \eqref{eq:tfy_prop}. The proof is divided into three parts wherein each part derives \eqref{eq:sprocD1}, \eqref{eq:sprocD2} and \eqref{eq:sprocD3}, respectively.
% If there exist matrices $U_\mathrm{e}$, $\bar{D}_1$, $\bar{D}_2$ and $\bar{D}_3$, with $\bar{D}_1+\bar{D}_2+\bar{D}_3 \succeq 0$ (cf. \eqref{eq:thm12}), and scalars $\tau_1 \geq 0$, $\tau_2 \geq 0$, and $\tau_3 \geq 0$ that satisfy \eqref{eq:sprocD1}, \eqref{eq:sprocD2}, and \eqref{eq:sprocD3}, the condition \eqref{eq:prop3} in Proposition \ref{prop4} is ensured.

\textbf{Part I.} Inequality \eqref{eq:D1} can be written as 
\begin{align}\label{eq:Vxtr_D1}
\nonumber
    \begin{bmatrix}
        (V_\mathrm{x,tr}^\mathsf{H} \otimes I_{n_\phi})\hat{Z}^\mathsf{H} \\ I
    \end{bmatrix}^\mathsf{H}&\\
\nonumber
    \times \begin{bmatrix}
        0 & \star^\mathsf{H} & \\
        \underbrace{\begin{bmatrix}
        (1-\epsilon)D_\mathrm{des}^{\frac{1}{2}\top}\left( \mathbf{1}_L^\top U_\mathrm{e}^\top \otimes I_{n_\phi} \right)\\
        0
    \end{bmatrix}}_{\star} & -\bar{D}_1
    \end{bmatrix}& \\
    \times \begin{bmatrix}
        (V_\mathrm{x,tr}^\mathsf{H} \otimes I_{n_\phi})\hat{Z}^\mathsf{H} \\ I
    \end{bmatrix} & \succeq 0.
\end{align}
% From \eqref{eq:tf_prop}, we have
% \begin{align*}
%     L\left( \left( \tilde{V}_\mathrm{x} \tilde{V}_\mathrm{x}^\mathsf{H}\right)\otimes  I_{n_\mathrm{x}n_\phi}\right)L^\mathsf{H} \preceq L \left( \tilde{\tilde{\Gamma}}_\mathrm{x} \otimes  I_{n_\mathrm{x}n_\phi} \right)L^\mathsf{H},
% \end{align*}
% \begin{align*}
%     V_\mathrm{\mathrm{x},tr} V_\mathrm{\mathrm{x},tr}^\mathsf{H}-V_\mathrm{\mathrm{x},tr} \hat{V}_\mathrm{x}^\mathsf{H} - \hat{V}_\mathrm{x} V_\mathrm{\mathrm{x},tr}^\mathsf{H} + \hat{V}_\mathrm{x} \hat{V}_\mathrm{x}^\mathsf{H} \preceq \tilde{\Gamma}_\mathrm{x}
% \end{align*} 
% and hence,
From \eqref{eq:vtilde} and \eqref{eq:tf_prop}, we have
\begin{align}
\nonumber
    & \hat{Z} \left(\left(  V_\mathrm{\mathrm{x},tr} V_\mathrm{\mathrm{x},tr}^\mathsf{H}-V_\mathrm{\mathrm{x},tr} \hat{V}_\mathrm{x}^\mathsf{H} - \hat{V}_\mathrm{x} V_\mathrm{\mathrm{x},tr}^\mathsf{H} \right)  \otimes  I_{n_\phi}\right) \hat{Z}^\mathsf{H}\\
    \preceq &\; \hat{Z} \left( \left(\tilde{\Gamma}_\mathrm{x}- \hat{V}_\mathrm{x} \hat{V}_\mathrm{x}^\mathsf{H}\right) \otimes  I_{n_\phi} \right)\hat{Z}^\mathsf{H},
\end{align}    
which can be equivalently written as
\begin{align}\label{eq:Vxtr_bound}
 \nonumber
    \begin{bmatrix}
        (V_\mathrm{x,tr}^\mathsf{H} \otimes I_{n_\phi})\hat{Z}^\mathsf{H} \\ I
    \end{bmatrix}^\mathsf{H}&\\
    \times \begin{bmatrix}
        -I & (\hat{V}_\mathrm{x}^\mathsf{H}  \otimes I_{n_\phi})\hat{Z}^\mathsf{H}\\
     \hat{Z}(\hat{V}_\mathrm{x} \otimes I_{n_\phi}) & \hat{Z}(( \tilde{\Gamma}_\mathrm{x}-\hat{V}_\mathrm{x} \hat{V}_\mathrm{x}^\mathsf{H})\otimes I_{n_\phi})\hat{Z}^\mathsf{H}
    \end{bmatrix}& \\
    \nonumber
    \times \begin{bmatrix}
        (V_\mathrm{x,tr}^\mathsf{H} \otimes I_{n_\phi})\hat{Z}^\mathsf{H} \\ I
    \end{bmatrix}&\succeq 0.
\end{align}   
By using the matrix S-lemma \cite{boyd2004convex,vanwaarde2022noisy}, Inequality \eqref{eq:Vxtr_D1} holds for all $V_\mathrm{x,tr}$ satisfying Inequality \eqref{eq:Vxtr_bound}, if $S_\textnormal{exp-1}(\epsilon, \tau_1, U_\mathrm{e}, \hat{Z}, \hat{V}_\mathrm{x}, \tilde{\Gamma}_\mathrm{x}, \tilde{D}_\mathrm{des}, \bar{D}_1)~\succeq~0$ \eqref{eq:sprocD1} holds with $\tau_1 \geq 0$.

\textbf{Part II.} Inequality \eqref{eq:D2} can be written as
\begin{align}\label{eq:Vphitr_D2}
\nonumber
    \begin{bmatrix}
        \begin{bmatrix}
        0\\ V_\mathrm{\phi,tr}\otimes I_{n_\mathrm{x} n_\phi}
    \end{bmatrix}^\mathsf{H} \\ I
    \end{bmatrix}^\mathsf{H}&\\
\nonumber
    \times \begin{bmatrix}
        0 & \underbrace{(1-\epsilon)(U_\mathrm{e}\otimes I_{n_\mathrm{x} n_\phi}) \hat{Z}^\mathsf{H}}_{\star}\\
        \star^\mathsf{H} & -(1-\epsilon)\hat{Z} \hat{Z}^\mathsf{H}-\left(\tfrac{1-\epsilon}{\epsilon}\right)\bar{W}_Z - \bar{D}_2
    \end{bmatrix}&\\
    \times \begin{bmatrix}
        \begin{bmatrix}
        0\\ V_\mathrm{\phi,tr}\otimes I_{n_\mathrm{x} n_\phi}
    \end{bmatrix}^\mathsf{H} \\ I
    \end{bmatrix} & \succeq 0.
\end{align}    

From \eqref{eq:vtilde} and \eqref{eq:tf_prop}, we have
% \begin{align}\label{eq:Vphitr}
%      \left( \tilde{V}_\mathrm{\phi}\tilde{V}_\mathrm{\phi}^\mathsf{H} \right)\otimes I_{n_\mathrm{x}^2n_\phi} - \tilde{\Gamma}_\phi \otimes I_{n_\mathrm{x}^2n_\phi} & \preceq 0.
% \end{align}
\begin{align}\label{eq:Vphitr}
\nonumber
    &(V_\mathrm{\phi,tr} V_\mathrm{\phi,tr}^\mathsf{H}-V_\mathrm{\phi,tr} \hat{V}_\phi^\mathsf{H} - \hat{V}_\phi V_\mathrm{\phi,tr}^\mathsf{H} + \hat{V}_\phi \hat{V}_\phi^\mathsf{H})\otimes I_{n_\mathrm{x} n_\phi}\\ \preceq & \; \tilde{\Gamma}_\phi \otimes I_{n_\mathrm{x} n_\phi},
\end{align} 

which can be equivalently written as
 \begin{align}\label{eq:Vphitr_bound}
 \nonumber
    \begin{bmatrix}
        \begin{bmatrix}
        0\\ V_\mathrm{\phi,tr}\otimes I_{n_\mathrm{x} n_\phi}
    \end{bmatrix}^\mathsf{H} \\ I
    \end{bmatrix}^\mathsf{H}&\\
    \times \begin{bmatrix}
        -I & \begin{bmatrix}
        0\\ \hat{V}_\mathrm{\phi}\otimes I_{n_\mathrm{x} n_\phi}
    \end{bmatrix}^\mathsf{H}\\
    \begin{bmatrix}
        0\\ \hat{V}_\mathrm{\phi}\otimes I_{n_\mathrm{x} n_\phi}
    \end{bmatrix} & \begin{bmatrix}
        0 & 0\\ 0 & (\tilde{\Gamma}_\phi - \hat{V}_\mathrm{\phi}\hat{V}_\mathrm{\phi}^\mathsf{H})\otimes I_{n_\mathrm{x} n_\phi}
    \end{bmatrix}
    \end{bmatrix}&\\
    \nonumber
    \times \begin{bmatrix}
        \begin{bmatrix}
        0\\ V_\mathrm{\phi,tr}\otimes I_{n_\mathrm{x} n_\phi}
    \end{bmatrix}^\mathsf{H} \\ I
    \end{bmatrix}&\succeq 0.
\end{align}   

By using the matrix S-lemma, Inequality \eqref{eq:Vphitr_D2} holds for all $V_\mathrm{\phi,tr}$ satisfying Inequality \eqref{eq:Vphitr_bound}, if $S_\textnormal{exp-2} (\epsilon,\tau_2, U_\mathrm{e},\hat{Z},\bar{W}_Z,\hat{V}_\phi, \tilde{\Gamma}_\phi,\tilde{D}_\mathrm{des},\bar{D}_2) \succeq 0$ \eqref{eq:sprocD2} holds with $\tau_2 \geq 0$.
% 

% \textbf{Part III.}
% Utilizing the structure of $\bar{D}_3$ \eqref{eq:slackDstructure}, Inequality \eqref{eq:D3} can be equivalently written as
% \begin{align}\label{eq:smallD3}
%         \left( (1-\epsilon) \tilde{\Phi} \tilde{\Phi}^\mathsf{H} - \left(\tfrac{1-\epsilon}{\epsilon}\right) \bar{W}_\phi\right) \otimes I_{n_\mathrm{x}} - \gamma_\mathrm{w} D_\mathrm{des} - \bar{D}_{3_{11}} & \succeq 0.
% \end{align}

% Inequality \eqref{eq:smallD3} can be written as
% \begin{scriptsize}
% \begin{align}\label{eq:Vphitr_D3}
% \nonumber
% \begin{bmatrix}
%          V_{\phi,\mathrm{tr}} \otimes I_{n_\mathrm{x}} \\ I
%     \end{bmatrix}^\mathsf{H} &\\
% \nonumber
%     \times \begin{bmatrix}
%         (1-\epsilon)(U_\mathrm{e} U_\mathrm{e}^\top)\otimes I_{n_\mathrm{x}} & 0\\
%         0 & -\left(\tfrac{1-\epsilon}{\epsilon}\right)\bar{W}_\phi-\gamma_\mathrm{w} D_\mathrm{des} -\bar{D}_{3_{11}} 
%     \end{bmatrix}&\\
%     \times \begin{bmatrix}
%          V_{\phi,\mathrm{tr}} \otimes I_{n_\mathrm{x}} \\ I
%     \end{bmatrix}& \succeq 0.    
% \end{align}
% \end{scriptsize}
% If $\bar{D}_3$ did not have the specified structure, 
\textbf{Part III.} Inequality \eqref{eq:D3} can be written as
\begin{align}\label{eq:Vphitr_D3_nostruc}
\nonumber
&\begin{bmatrix}
    \begin{bmatrix}
            V_{\phi,\mathrm{tr}}\\0
        \end{bmatrix}^\mathsf{H} \\ I
    \end{bmatrix}^\mathsf{H} &\\
\nonumber
    & \times \begin{bmatrix}
        (1-\epsilon)U_\mathrm{e} U_\mathrm{e}^T & 0\\
        0 & \begin{bmatrix}
            -\left(\tfrac{1-\epsilon}{\epsilon}\right)\bar{W}_\phi-\frac{\gamma_\mathrm{w}}{T} D_\mathrm{des}  & 0\\
            0 & 0
        \end{bmatrix} -\bar{D}_{3} 
    \end{bmatrix}&\\
\nonumber
    & \times \begin{bmatrix}
        \begin{bmatrix}
            V_{\phi,\mathrm{tr}}\\0
        \end{bmatrix}^\mathsf{H} \\ I
    \end{bmatrix}\\
    \succeq & \; 0.    
\end{align}

% \begin{align}\label{eq:VphitrInx}
%      \left( \tilde{V}_\mathrm{\phi}\tilde{V}_\mathrm{\phi}^\mathsf{H} \right)\otimes I_{n_\mathrm{x}} - \tilde{\Gamma}_\phi \otimes I_{n_\mathrm{x}} & \preceq 0.
% \end{align}
From \eqref{eq:vtilde} and \eqref{eq:tf_prop}, we have
\begin{align}
V_\mathrm{\phi,tr} V_\mathrm{\phi,tr}^\mathsf{H}-V_\mathrm{\phi,tr} \hat{V}_\phi^\mathsf{H} - \hat{V}_\phi V_\mathrm{\phi,tr}^\mathsf{H} + \hat{V}_\phi \hat{V}_\phi^\mathsf{H}\preceq \tilde{\Gamma}_\phi,
\end{align} 
which can be equivalently written as
% \begin{align}\label{eq:VphitrInx_bound}
%  \nonumber
%     \begin{bmatrix}
%          V_\mathrm{\phi,tr}^\mathsf{H} \otimes I_{n_\mathrm{x}} \\ I
%     \end{bmatrix}^\mathsf{H}&\\
% \nonumber
%     \times \begin{bmatrix}
%         -I & \hat{V}_\mathrm{\phi}^\mathsf{H} \otimes I_{n_\mathrm{x}} \\
%         \hat{V}_\mathrm{\phi}\otimes I_{n_\mathrm{x}}
%      &  (\tilde{\Gamma}_\phi - \hat{V}_\mathrm{\phi}\hat{V}_\mathrm{\phi}^\mathsf{H}) \otimes I_{n_\mathrm{x}}
%     \end{bmatrix}\\
%     \times \begin{bmatrix}
%         V_\mathrm{\phi,tr}^\mathsf{H} \otimes I_{n_\mathrm{x}} \\ I
%     \end{bmatrix}&\succeq 0.
% \end{align}   
% or
\begin{align}\label{eq:VphitrInx_bound_kron}
 \nonumber
    \begin{bmatrix}
         \begin{bmatrix}
             V_\mathrm{\phi,tr} \\0
         \end{bmatrix}^\mathsf{H} \\ I
    \end{bmatrix}^\mathsf{H}
    \begin{bmatrix}
        -I & \begin{bmatrix}
            \hat{V}_\mathrm{\phi}\\0
        \end{bmatrix}^\mathsf{H}  \\
        \begin{bmatrix}
            \hat{V}_\mathrm{\phi}\\0
        \end{bmatrix}
     &  \begin{bmatrix}(\tilde{\Gamma}_\phi - \hat{V}_\mathrm{\phi}\hat{V}_\mathrm{\phi}^\mathsf{H}) & 0 \\ 0 & 0 \end{bmatrix}
    \end{bmatrix}\\
    \times \begin{bmatrix}
       \begin{bmatrix}
             V_\mathrm{\phi,tr} \\0
         \end{bmatrix}^\mathsf{H} \\ I
    \end{bmatrix}&\succeq 0.
\end{align} 

% By using the matrix S-lemma, Inequality \eqref{eq:Vphitr_D3} holds for all $V_\mathrm{\phi,tr}$ satisfying Inequality \eqref{eq:VphitrInx_bound}, if and only if the following holds with $\tau_3 \geq 0$:
% \begin{footnotesize}
% \begin{align}\label{eq:presproc_D3}
% \nonumber
%     \begin{bmatrix}
%         (1-\epsilon)(U_\mathrm{e} U_\mathrm{e}^\top)\otimes I_{n_\mathrm{x}} & 0\\
%         0 & -\left(\tfrac{1-\epsilon}{\epsilon}\right)\bar{W}_\phi-\gamma_\mathrm{w} D_\mathrm{des} -\bar{D}_{3_{11}}
%     \end{bmatrix}&\\
%     -\tau_3 \begin{bmatrix}
%          -I & \hat{V}_\mathrm{\phi}^\mathsf{H} \otimes I_{n_\mathrm{x}} \\
%         \hat{V}_\mathrm{\phi}\otimes I_{n_\mathrm{x}} &  (\Gamma_\phi - \hat{V}_\mathrm{\phi} \hat{V}_\mathrm{\phi}^\mathsf{H}) \otimes I_{n_\mathrm{x}}
%     \end{bmatrix}& \succeq 0.
% \end{align}   
% \end{footnotesize}

By using the matrix S-lemma, Inequality \eqref{eq:Vphitr_D3_nostruc} holds for all $V_\mathrm{\phi,tr}$ satisfying Inequality \eqref{eq:VphitrInx_bound_kron}, if the following holds with $\tau_3 \geq 0$:
\begin{align}\label{eq:presproc_D3_nostruc}
\nonumber
    \begin{bmatrix}
        (1-\epsilon) U_\mathrm{e} U_\mathrm{e}^T & 0\\
        0 & \begin{bmatrix}
            -\left(\tfrac{1-\epsilon}{\epsilon}\right)\bar{W}_\phi-\frac{\gamma_\mathrm{w}}{T} D_\mathrm{des}  & 0\\
            0 & 0
        \end{bmatrix} -\bar{D}_{3} 
    \end{bmatrix}&\\
    -\tau_3 \begin{bmatrix}
        -I & \begin{bmatrix}
            \hat{V}_\mathrm{\phi}\\0
        \end{bmatrix}^\mathsf{H}  \\
        \
        \begin{bmatrix}
            \hat{V}_\mathrm{\phi}\\0
        \end{bmatrix}
     &  \begin{bmatrix}(\tilde{\Gamma}_\phi - \hat{V}_\mathrm{\phi}\hat{V}_\mathrm{\phi}^\mathsf{H}) & 0 \\ 0 & 0 \end{bmatrix}
    \end{bmatrix} & \succeq 0.
\end{align}

For matrices $U_\mathrm{e}$ and $\hat{U}$, we have $\bar{U}_\mathrm{e}\bar{U}_\mathrm{e}^\mathsf{H} - \bar{U}_\mathrm{e}\hat{U}^\mathsf{H}- \hat{U} \bar{U}_\mathrm{e}^\mathsf{H} + \hat{U} \hat{U}^\mathsf{H} = (\bar{U}_\mathrm{e} - \hat{U})(\bar{U}_\mathrm{e} - \hat{U})^\mathsf{H}\succeq 0$, and hence, we have
\begin{align}\label{eq:uuconvex}
    U_\mathrm{e} U_\mathrm{e}^T \succeq  U_\mathrm{e} \hat{U}^\top + \hat{U} U_\mathrm{e}^\top - \hat{U} \hat{U}^\top.
\end{align}

Inserting Inequality \eqref{eq:uuconvex} in Inequality \eqref{eq:presproc_D3_nostruc} yields $S_\textnormal{exp-3}(\epsilon, \tau_3, U_\mathrm{e},\hat{U},\hat{W}_\phi,\hat{V}_\phi, \tilde{\Gamma}_\phi,D_\mathrm{des},\bar{D}_3) \succeq 0$ \eqref{eq:sprocD3}.

Therefore, if there exist matrices $U_\mathrm{e}$, $\bar{D}_1$, $\bar{D}_2$ and $\bar{D}_3$ that satisfy Inequalities \eqref{eq:sprocD1}, \eqref{eq:sprocD2}, \eqref{eq:sprocD3} with $\bar{D}_3$, with $\bar{D}_1+\bar{D}_2+\bar{D}_3 \succeq 0$ (cf. \eqref{eq:thm12}), then Inequality \eqref{eq:prop3} is satisfied for all $V_{\mathrm{x,tr}}$, $V_{\mathrm{\phi,tr}}$, $Y_{\mathrm{x,tr}}$, and $Y_{\mathrm{\phi,tr}}$ satisfying bounds \eqref{eq:tf_prop}, \eqref{eq:tfy_prop}, and the exploration goal \eqref{eq:exp_goal2} is achieved. $\hfill\square$
\end{pf}}{}

% DONT UNCOMMENT AFTER HERE
% \section{Sample-based constants $\tilde{\Gamma}_\mathrm{x}$, $\Gamma_\mathrm{x}$ and $\Gamma_\phi$ }\label{app:scenario}
% \input{sections/scenario}

% \section{Robust constants $\tilde{\Gamma}_\mathrm{x}$, $\tilde{\Gamma}_\phi$, $\Gamma_\mathrm{x}$ and $\Gamma_\phi$ }\label{app:robustbounds}
% \input{sections/robustbounds}

\end{document}